\documentclass[12pt,reqno]{amsart}

\usepackage{fullpage,amsmath,amssymb,amsfonts,mathrsfs,bm,graphicx,hyperref}
\usepackage{enumitem}

\def\bbone{{\mathchoice {\rm 1\mskip-4mu l} {\rm 1\mskip-4mu l}
{\rm 1\mskip-4.5mu l} {\rm 1\mskip-5mu l}}}

\newcommand{\ua}{\mathbf{a}}
\newcommand{\ux}{\mathbf{x}}
\newcommand{\vepi}{{\varepsilon_{i \ast}}}
\newcommand{\taui}{\tau_{i \ast}}
\newcommand{\wtA}{\widetilde{A}}
\newcommand{\wtM}{\widetilde{M}}
\newcommand{\wts}{\widetilde{\sigma}}
\newcommand{\wtmu}{\widetilde{\mu}}
\newcommand{\wtla}{\widetilde{\lambda}}

\newtheorem{theorem}{Theorem}[section]
\newtheorem{remark}{Remark}[section]

\newtheorem{corollary}{Corollary}[section]
\newtheorem{lemma}{Lemma}[section]
\newtheorem{proposition}{Proposition}[section]

\begin{document}

\author{Abdelmalek Abdesselam}

\address{Abdelmalek Abdesselam, Department of Mathematics,
P. O. Box 400137,
University of Virginia,
Charlottesville, VA 22904-4137, USA}

\email{malek@virginia.edu}

\title{A combinatorial formula for the coefficients of multidimensional resultants}

\begin{abstract}
The classical multidimensional resultant can be defined as the, suitably normalized, generator of a projective elimination ideal in the ring of universal coefficients. This is the approach via the so-called inertia forms or Tr\"{a}gheitsformen. Using clever substitutions, Mertens and Hurwitz gave a criterion, for recognizing such inertia forms, which amounts to a linear system for their numerical coefficients. In this article we explicitly solve this linear system. We do so by identifying a subset of the available equations which forms a unitriangular system. The key notion we use is that of transversal, i.e., a selection of a monomial term in each of the homogeneous polynomials at hand. We need two such transversals which are disjoint and extremal, in the sense that they relate to extremizers of a, possibly new, determinantal inequality for differences of two substochastic matrices. Thanks to this notion of extremal pair of transversals, we derive an explicit formula for the coefficients of general multidimensional resultants, as a sum of terms made of a sign times a product of multinomial coefficients, thereby explicitly showing they are integer-valued.
As an application of our formula, we recover Sombra's bound on the height of resultants, in the classical case.
\end{abstract}

\maketitle

\tableofcontents

\section{Introduction}

For $n\ge 1$, and temporarily working over the field $\mathbb{C}$ of complex numbers, consider $n$ homogeneous polynomials (or {\it forms})
$F_1(x),\ldots,F_n(x)$, in $n$ variables $x_1,\ldots,x_n$, of respective degrees $d_1,\ldots,d_n\ge 1$.
There is a unique polynomial ${\rm Res}_{d_1,\ldots,d_n}(F_1,\ldots,F_n)$ in the coefficients of $F_1,\ldots,F_n$ which satisfies the following axioms.
\begin{itemize}
\item
{\bf A1:}
${\rm Res}_{d_1,\ldots,d_n}(F_1,\ldots,F_n)=0$ iff $\exists x\in\mathbb{C}^n\backslash\{0\}, \forall i, F_i(x)=0$.
\item
{\bf A2:}
$\forall i$, ${\rm Res}_{d_1,\ldots,d_n}(F_1,\ldots,F_n)$ is homogeneous of degree $\prod_{j\neq i}d_j$ in the coefficients of $F_i$.
\item
{\bf A3:}
${\rm Res}_{d_1,\ldots,d_n}(F_1,\ldots,F_n)=1$ when $F_1(x)=x_1^{d_1},\ldots,F_n(x)=x_n^{d_n}$ (the `diagonal case').
\end{itemize}
It turns out, ${\rm Res}_{d_1,\ldots,d_n}$ has {\it integer coefficients}. This allows one to feed it the coefficients of homogeneous polynomials over any field $K$, of {\it arbitrary characteristic}. Let $\overline{K}$ denote the algebraic closure of $K$, then the analogue of axiom {\bf A1} holds for forms with coefficients in $K$: ${\rm Res}_{d_1,\ldots,d_n}(F_1,\ldots,F_n)=0$ iff $\exists x\in{\overline{K}}^n\backslash\{0\}, \forall i, F_i(x)=0$.
The polynomial ${\rm Res}_{d_1,\ldots,d_n}$ therefore is a universal object, of fundamental importance in the study of polynomial systems of equations, and it is 
called the {\it multidimensional resultant}.
If $(d_1,\ldots,d_n)=(1,\ldots,1)$, then the resultant is just the familiar determinant of the matrix formed by the coefficients of the linear forms $F_1,\ldots,F_n$.
For modern introductions to resultants, their construction, and their properties, the reader may consult~\cite{CattaniD,SturmfelsCBMS,CoxLOusing,ElkadiM,Jouanolou,GKZbook}.
More elementary or classical presentations can be found in~\cite{FaadB,Macaulaybook,VdWaerden,HodgeP}. Although the study of resultants goes back several centuries, and although they are polynomials, to the best of our knowledge, no explicit polynomial formula, for general $n$ and degree sequence $d_1,\ldots,d_n$, was available until the recent work~\cite{MorozovS}.
Prior to that, explicit polynomial formulas (e.g., determinantal) were only known in particular cases (see the list in~\cite[Ch. 13, Prop. 1.6]{GKZbook}, later improved in~\cite{EisenbudSW}).
Formulas which were valid for any $n$ and degree sequence, were typically given as a ratio of products of determinants, i.e., the determinant of a complex (see~\cite[Ch. 13, Thm. 1.4]{GKZbook}), or following Macaulay~\cite{MacaulayFormulas}, as a ratio of two large determinants (see~\cite{DAndreaD} for recent generalizations).

In this article, we obtain a new explicit polynomial formula for the resultant ${\rm Res}_{d_1,\ldots,d_n}$ that is valid for arbitrary $n$ and degree sequence. Moreover, we explicitly obtain all the coefficients of these resultants as alternating sums of products of multinomial coefficients. It is thus manifest from our formula that multidimensional resultants have integer coefficients. The earlier Morozov-Shakirov formula~\cite[\S7]{MorozovS} gives ${\rm Res}_{d_1,\ldots,d_n}$ as a polynomial in expressions called ``traces'', which themselves are polynomials in the coefficients of the forms $F_i$. It is a nonlinear generalization of the formula
\[
\left[ \exp\left({\rm tr}\log(I+A)\right)
\right]_n
\]
for the determinant of an $n\times n$ matrix $A$ of formal variables, where $[\cdots]_n$ means taking the homogeneous part of degree $n$, and the exponential and logarithm are understood as formal power series. A formula for the individual coefficients of the resultant is not given in~\cite{MorozovS}, although one should, in principle, be able to do so. However, since there are products of factorials in the denominators, it is not clear from the Morozov-Shakirov formula, that the coefficients of the resultant are integers, in contrast to our new formula.

In the following, we will give an outline of our approach, while referring to \S\ref{notationsec}
for our precise choice of notation.
Our starting point is an alternative definition of resultants via the notion of {\it inertia forms} (or Tr\"{a}gheitsformen). We write our, now generic, homogeneous polynomials $F_i$, $1\le i\le m$,  as
\[
F_i=\sum_{|\alpha|=d_i} a_{i,\alpha}\ \ux^{\alpha}\ ,
\]
where the $a_{i,\alpha}$ are {\it formal variables}, collectively denoted by $\ua=(a_{i,\alpha})_{(i,\alpha)\in\mathbb{D}}$.
Likewise, the $x_i$, $1\le i\le n$, are now considered as formal variables and are collectively denoted by $\ux$. Note that we, temporarily, allow a number of forms $m$ which may be different from $n$. 
We will primarily work in the polynomial ring $\mathbb{Q}[\ua,\ux]$. The set $\mathfrak{I}$ of inertia forms is that of elements $\mathscr{R}\in\mathbb{Q}[\ua,\ux]$ for which there exists $q\ge 0$, such that for all multiindex $\gamma$ with $|\gamma|=q$, there exist polynomials $G_1,\ldots,G_m$ in $\mathbb{Q}[\ua,\ux]$, such that
\begin{equation}
\ux^{\gamma}\ \mathscr{R}=F_1\ G_1+\cdots+F_m\ G_m\ .
\label{Bezout}
\end{equation}
This is an ideal in the ring $\mathbb{Q}[\ua,\ux]$, but we will focus on $\mathfrak{I}_0:=\mathfrak{I}\cap \mathbb{Q}[\ua]$, which is an ideal in the subring $\mathbb{Q}[\ua]$. An element of $\mathfrak{I}_0$ is called an {\it inertia form of degree zero} (IFDZ), where the degree here refers to the $x$ variables. Thus
\[
\mathfrak{I}_0=\left(\ 
\langle F_1,\ldots,F_m\rangle \ :\ \langle x_1,\ldots,x_n\rangle^{\infty}
\ \right)\cap \mathbb{Q}[\ua]
\]
is a particular case of projective elimination ideal.
However, it is a very special particular case because of the availability of what we call {\it Mertens-Hurwitz (MH) substitutions}. Suppose that for each $i$, we choose a monomial $\ux^\vepi$ of degree $d_i$ featuring in the form $F_i$ which we write as
\[
F_i=a_{i,\vepi}\ \ux^{\vepi}+\widetilde{F}_i\ ,
\]
where $\widetilde{F}_i$ denotes the sum of all the other terms in $F_i$. An MH substitution is the result of setting, for all $i$, $1\le i\le m$,
\begin{equation}
a_{i,\vepi}:=-\ \frac{\widetilde{F}_i}{\ux^{\vepi}}\ .
\label{MHsub}
\end{equation}
The set $\mathbb{E}=\{(i,\vepi)\ |\ 1\le i\le m\}$ which encodes the monomial selection is what we call a {\it transversal} and it is a subset of the set $\mathbb{D}$, or {\it diagram} which labels the $a_{i,\alpha}$ variables.
This substitution gives a ring homomorphism
\[
\mathbb{Q}[\ua,\ux]=\mathbb{Q}[\ua_{\mathbb{D}},\ux]
\longrightarrow 
\mathbb{Q}[\ux,\ux^{-1}][\ua_{\mathbb{D}\backslash\mathbb{E}}]\ ,
\]
into the ring of polynomials in the remaining $a_{i,\alpha}$ variables over the ring of Laurent polynomials in the $x$ variables, with rational coefficients.
The result of this substitution on the polynomial $\mathscr{R}\in\mathbb{Q}[\ua]$ featuring in the B\'ezout relation (\ref{Bezout})
will be denoted by $\overline{\mathscr{R}}^{\mathbb{E}}$. 
Multiplying $\overline{\mathscr{R}}^{\mathbb{E}}$ by a suitable power of the $x$'s in order to clear denominators, we then get a polynomial $\widehat{\mathscr{R}}^{\mathbb{E}}\in \mathbb{Q}[\ua_{\mathbb{D}\backslash\mathbb{E}},\ux]$.
We now have that, for a polynomial $\mathscr{R}\in\mathbb{Q}[\ua]$, being an IFDZ is
equivalent to $\widehat{\mathscr{R}}^{\mathbb{E}}$ being identically zero.
Moreover, this equivalence holds regardless of the prior choice of transversal $\mathbb{E}$. An equivalent definition of the resultant ${\rm Res}_{d_1,\ldots,d_n}$ is that it is the unique IFDZ, in the $m=n$ situation, which is multihomogeneous of multidegree $(\delta_1,\ldots,\delta_n)$ with $\delta_i=\prod_{j\neq i}d_j$, and which is equal to 1 when one plugs into it the coefficients of the diagonal collection of forms $F_i=x_i^{d_i}$.
Writing the expansion of this IFDZ as
\[
\mathscr{R}=\sum_{A} r_A\ \ua^{A}\ ,
\]
we see that finding the resultant amounts to solving a large system of linear equations
in the unknowns $r_A$, i.e., a problem of linear algebra over $\mathbb{Q}$. The main innovation in this article is a {\it method} for explicitly solving such a system.

Inertia forms were introduced in the work of Mertens and Hurwitz~\cite{Mertens1,Mertens2,Mertens3,Hurwitz}. For an introduction to their theory and the equivalence with the previous definition of resultants with axioms {\bf A1}-{\bf A3}, see~\cite[Ch. XI]{VdWaerden} and~\cite[Ch. IV]{HodgeP} from a `semiclassical' standpoint, and~\cite{JouanolouId,Jouanolou} from a modern perspective.
See also~\cite{AlAmrani} for a pedagogical introduction to the last two references by Jouanolou. The reader may also find~\cite[Ch. 6]{Buse}, and~\cite[Ch. 7]{Ferretti} quite helpful, as first introductions to inertia forms.
Restricting to the situation $m=n$ (the principal case of elimination theory), the system to be solved is made of equations denoted by ${\rm Eq}(\mathbb{E},C,\beta)$ which say:
\begin{equation}
\left[\ua^{C}\ux^{\beta}\right]\ \widehat{\mathscr{R}}^{\mathbb{E}}=0\ ,
\label{shortmaineq}
\end{equation}
namely, express the vanishing of the coefficient of the monomial $\ua^{C}\ux^{\beta}$ in the polynomial $\widehat{\mathscr{R}}^{\mathbb{E}}\in \mathbb{Q}[\ua_{\mathbb{D}\backslash\mathbb{E}},\ux]\subset\mathbb{Q}[\ua_{\mathbb{D}},\ux]$. Our approach is to identify, among the large number of equations ${\rm Eq}(\mathbb{E},C,\beta)$, a subset which essentially forms a {\it unitriangular system} which can readily be solved by back substitution. For this, we need two tightly entangled ingredients. We need a partial order among the $A$'s (or fillings) which store the exponents
in the $a_{i,\alpha}$ variables, so that $\widetilde{A}\prec A$ implies we must solve for the unknown $r_{\widetilde{A}}$ before $r_A$. We also need, for each $A$, a prescription for choosing an associated {\it resolvent equation} ${\rm Eq}(\mathbb{E},C,\beta)$ of the form
\begin{equation}
\lambda\ r_A+\sum_{\widetilde{A}\prec A} \rho_{\widetilde{A}}\ r_{\widetilde{A}}=0\ .\label{lambdarhoeq}
\end{equation}
This is the equation we will use in order to solve for $r_{A}$, in terms of $r_{\widetilde{A}}$'s which, at this point, should already be known. Of course, it is essential to avoid vicious circles in this process.
Our only guide for choosing these resolvent equations is sheer optimism. Since we do not want to introduce denominators, so the solutions $r_A$ can be shown to be integers, we look for equations where $\lambda=\pm 1$. This amounts to appropriately choosing {\it another transversal} $\mathbb{T}$ which is disjoint from $\mathbb{E}$. While the proper choice of transversal pair $(\mathbb{E},\mathbb{T})$ could a priori depend on $A$ when deciding how to solve for the unknown $r_A$, we found, by inspection, a fixed pair which works. The transversal $\mathbb{E}$ corresponds to the matrix
\begin{equation}
\varepsilon=\begin{pmatrix}
d_1 & 0 & \cdots & 0 \\
0 & d_2 & \ddots & \vdots \\
\vdots & \ddots & \ddots & 0 \\
0 & \cdots & 0 & d_n
\end{pmatrix}
\label{diagonaleps}
\end{equation}
where $\vepi=(\varepsilon_{ij})_{1\le j\le n}$ is stored as the $i$-th row of the above matrix.
This means that we select the monomials $x_1^{d_1},x_2^{d_2},\ldots,x_n^{d_n}$ for $i=1,2,\ldots,n$ respectively, when performing the MH substitutions.
Likewise, the transversal $\mathbb{T}$ corresponds to the matrix
\begin{equation}
\tau=\begin{pmatrix}
d_1-1 & 1 & 0 & \cdots & 0 \\
d_2-1 & 0 & 1 & \ddots & \vdots \\
\vdots & \vdots & \ddots & \ddots & 0 \\
d_{n-1}-1 & 0 & \cdots & 0 & 1 \\
d_n & 0 & \cdots & 0 & 0
\end{pmatrix}
\label{taupick}
\end{equation}
whose rows we will denote by $\taui:=(\tau_{ij})_{1\le j\le n}$, and which can also be seen as the monomial selection
\[
x_1^{d_1-1}x_2\ ,\  x_1^{d_2-1}x_3\ ,\ \ldots\ ,\  x_{1}^{d_{n-1}-1} x_{n}\ ,\  x_1^{d_n}\ ,
\]
for $i=1,2,\ldots,n$ respectively. 
The above pair $(\mathbb{E},\mathbb{T})$ is an example of what we call an {\it extremal pair of transversals} which is the key concept in this article. The corresponding matrices $\varepsilon,\tau$ are related to extremizers which saturate the following determinantal inequality which may be new. If $P,Q$ are two square substochastic matrices, then
\begin{equation}
|{\rm det}(P-Q)|\le 1\ .
\label{substoch}
\end{equation}
With the right choice of resolvent equation ${\rm Eq}(\mathbb{E},C,\beta)$, we can write
\[
r_A=\sum_{\widetilde{A}\prec A}\mathscr{T}(A,\widetilde{A})\ r_{\widetilde{A}}\ ,
\]
for a suitable {\it transition matrix} $\mathscr{T}$.
Similarly to the computation of the probability of hitting a specific absorbing state, in elementary Markov chain theory, we iterate the above equation, until we hit fillings $A$ which are {\it first-row-reduced}, namely have all $A_{1,\alpha}=0$, except when $\alpha=(d_1,0,\ldots,0)$.
At this point we suspend the linear system resolution, and use a nonlinear shortcut.
The product formula~\cite[\S5.7]{Jouanolou} and Laplace formula~\cite[\S5.10]{Jouanolou} (see also the reduction by one variable in~\cite[\S6.4]{Buse}) for the resultant produce the identity
\begin{equation}
{\rm Res}_{d_1,\ldots,d_n}\left(x_1^{d_1},\overline{F}^{(2)}_2,\ldots,\overline{F}^{(2)}_n\right)
=
\left({\rm Res}_{d_2,\ldots,d_n}\left(\overline{F}^{(2)}_2,\ldots,\overline{F}^{(2)}_n\right)\right)^{d_1}\ ,
\label{Laplaceeq}
\end{equation}
where the $\overline{F}^{(2)}_i$ are obtained from the $F_i$ by removing all the monomials which involve the first variable $x_1$.
From this identity, one easily obtains, by the multinomial theorem, an expansion for $r_A$, when $A$ is first-row-reduced, in terms of products of $d_1$ quantities like $r_A$ but {\it in one dimension less}. A straightforward iteration of the above procedure, by induction on the dimension $n$, gives $r_A$ as an explicit combinatorial sum, along a tree of fixed shape, of terms made of a sign times a product of multinomial coefficients. The tree has a root with $d_1$ children, which each have $d_2$ children, etc. This explicit formula for $r_A$ is given by Theorem 
\ref{maintheorem} in \S\ref{formulasec} below, and it is the main new result of the present article.

In \S\ref{notationsec}, we introduce the notations and definitions which we will use throughout this article. We will also review some basic properties of inertia forms.
In \S\ref{warmupsec}, we explain what the partial order $\prec$ is about, and use our transversal pair in order to show, when $m=n$, that if the degree in one of the forms is less than $\prod_{j\neq i}d_j$, then there is no nonzero IFDZ of the corresponding multidegree. A consequence is that if $m<n$, there are no nonzero IFDZ's. These are known results, but the point is to see how our general method works in a simpler situation than the computation of the resultant itself. We also prove, when $m=n$, using our new method, that in the correct multidegree, a nonzero IFDZ must have balanced weight, as expected from ${\rm SL}_n$-invariance. 
In \S\ref{convexsec}, we take a short break from heavy explicit commutative algebra, and examine the notion of extremal pair of transversals, and prove the determinantal inequality (\ref{substoch}). In \S\ref{formulasec}, we derive our explicit formula for the $r_A$'s.
We then look at the detailed implementation of our general formula in the examples provided by the linear case in \S\ref{linearsec}, and the case of two binary forms in \S\ref{binarysec}. In \S\ref{heightsec}, we give a new proof of a bound due to Sombra for the height of resultants, in the classical case. In \S\ref{outlooksec}, we finally conclude our article with a brief discussion of outlook and open problems left for future work.

\smallskip\noindent
\textit{\textbf{Acknowledgements:}} We thank Carlos D'Andrea, Darij Grinberg, and Eric Rowland for valuable feedback on the first draft of this article. We also thank the anonymous referees for helpful comments and suggestions.

\section{Notations, definitions, and basic facts about inertia forms}
\label{notationsec}

The object of this section is to introduce the notation needed for the rest of this article. We first use it in \S\ref{notationsec} and \S\ref{warmupsec} in order to recover some well known results about inertia forms, going back to Hurwitz. We caution the reader that the notation introduced in \S\ref{notationsec} is much more complicated than needed for the basic results reviewed in \S\ref{notationsec} and \S\ref{warmupsec}, which can be established by simpler means. However, our notation which may prima facie seem excessively heavy is essentially unavoidable when formulating and deriving our main result, i.e., Theorem \ref{maintheorem}. To arrive at the explicit formulas for the resultant coefficients $r_A$ in that theorem, we must first explicitly write the linear system of equations which the $r_A$ satisfy. The individual equations of that system, denoted by ${\rm Eq}(\mathbb{E},C,\beta)$, must therefore be stated explicitly, and this requires the heavy notation to be introduced in this section. The point of \S\ref{notationsec} and \S\ref{warmupsec} is to give the reader a chance to get accustomed to our notation and get some practice with it, on simpler problems where our notation is not necessary, before the more involved use made in \S\ref{formulasec} where our notation, or some equivalent, becomes necessary. We tried to make the notation easier to follow by using consistency in the choice of symbols and characters. For instance, exponents of monomials in the $x$ variables are encoded by multiindices typically denoted by 
Greek letters $\alpha,\beta,\gamma,\ldots$ The coefficients of the generic forms $F$, i.e., the variables $a_{i,\alpha}$ are indexed by a finite set $\mathbb{D}$ defined below. Subsets of the latter, like our transversals $\mathbb{E},\mathbb{T}$, are typically denoted by blackboard characters. Monomials in the $a_{i,\alpha}$ variables are encoded by more complicated multiindices we call fillings and typically denoted by capital Roman letters like $A,B,C,M$.

In this article we adopt the convention $\mathbb{N}:=\{0,1,2,\ldots\}$, and if we need the set of strictly positive integers, we will denote it by $\mathbb{Z}_{>0}$.
For any $q\in\mathbb{N}$, we use the standard notation $[q]:=\{1,2,\ldots,q\}$, so that $|[q]|=q$, where $|X|$ denotes the cardinality of a finite set $X$. The set of maps $f:X\rightarrow Y$, from a set $X$ to a set $Y$, will be denoted by $Y^{X}$. 
If we write an inclusion $X\subset Y$, we allow the sets to be equal.
A {\it format} is a triple $(m,n,d)$ where $m,n\in\mathbb{Z}_{>0}$ and $d=(d_1,\ldots,d_m)^{\rm T}$ is a {\it column} vector made of the integers $d_1,\ldots,d_m\in\mathbb{Z}_{>0}$ which will serve as the degrees of our generic forms of interest $F_1,\ldots,F_m$.
Given such a format, we define the associated diagram $\mathbb{D}$, as follows.
Our terminology is borrowed from the theory of Young diagrams and tableaux. We let
\[
\mathbb{D}:=\{(i,\alpha)\ |\ i\in[m]\ {\rm and}\ \alpha\in\mathbb{N}^n,\ |\alpha|=d_i\}\ .
\]
An element $(i,\alpha)\in\mathbb{D}$ will be called a {\it cell} of the diagram.
Its purpose is to label the formal variable $a_{i,\alpha}$, as in the introduction.
An element $\alpha=(\alpha_1,\ldots,\alpha_n)\in\mathbb{Z}^n$
will be called a {\it weight vector}, and $|\alpha|$ means $\alpha_1+\cdots+\alpha_n$.
In the particular case where $\alpha\in\mathbb{N}^n$, we will prefer to call $\alpha$ a {\it multiindex}.
Note that $|\cdot|$ will be rather overused notation in this article, but it should be clear from the context if we mean cardinality of a finite set, the lenght of a multiindex (or weight vector), or, later in \S\ref{binarysec}, the weight of an integer partition, or simply the absolute value/modulus of a number.

We call $\mathbb{D}_i:=\mathbb{D}\cap(\{i\}\times\mathbb{N}^n)$
the {\it $i$-th row} of the diagram $\mathbb{D}$. This is related to the set $\widetilde{\mathbb{D}}_i$ of multiindices $\alpha$ such that $|\alpha|=d_i$.
The latter set can also be thought of as the set of monomials $\ux^{\alpha}:=x_1^{\alpha_1}\cdots x_n^{\alpha_n}$ of degree $d_i$.
A map $A:\mathbb{D}\rightarrow\mathbb{Z}$, $(i,\alpha)\mapsto A_{i,\alpha}$,  will be called a {\it filling} of the diagram $\mathbb{D}$. If $A\in\mathbb{N}^{\mathbb{D}}\subset\mathbb{Z}^{\mathbb{D}}$, we call $A$ an {\it effective filling}.
The main use of an effective filling $A$ is to do precise computations with monomials
\[
\ua^A:=\prod_{(i,\alpha)\in\mathbb{D}}a_{i,\alpha}^{A_{i,\alpha}}\ ,
\]
in the formal variables packaged in $\ua=\ua_{\mathbb{D}}:=(a_{i,\alpha})_{(i,\alpha)\in\mathbb{D}}$.
For a filling $A$, we define the {\it row sum} column vector
\[
{\rm Row}(A):=\begin{pmatrix}
{\rm Row}_1(A) \\
\vdots \\
{\rm Row}_m(A)
\end{pmatrix}\ ,
\]
where $\forall i\in[m], {\rm Row}_i(A):=\sum_{\alpha\in\widetilde{\mathbb{D}}_i} A_{i,\alpha}$. We also use the notation
\[
\|A\|:=\sum_{(i,\alpha)\in\mathbb{D}}A_{i,\alpha}\ ,
\]
as relief for the overworked $|\cdot|$.
We also define the associated weight vector
\[
{\rm wt}(A):=\sum_{(i,\alpha)\in\mathbb{D}}A_{i,\alpha}\ \alpha\ \in\ \mathbb{Z}^n\ .
\]
Note that we have the obvious identity
\begin{equation}
|{\rm wt}(A)|=\sum_{(i,\alpha)\in\mathbb{D}}A_{i,\alpha}\ |\alpha|=\sum_{i=1}^{m}
{\rm Row}_i(A)\ d_i={\rm Row}(A)^{\rm T}d\ .
\label{obviouswteq}
\end{equation}
When computing with weight vectors, it will be convenient to use the standard basis vectors $e_j=(0,\ldots,0,1,0,\ldots,0)$, $1\le j\le n$, where the 1 is in the $j$-th position.
If $\rho=(\rho_1,\ldots,\rho_m)^{\rm T}$ is a column vector with nonnegative integer entries, and if $A$ is an effective filling such that ${\rm Row}(A)=\rho$, we define
\begin{equation}
\begin{bmatrix}
\rho \\
A
\end{bmatrix}
:=\frac{\prod_{i=1}^m \rho_i!}{\prod_{(i,\alpha)\in\mathbb{D}}A_{i,\alpha}!}
\label{supernomial}
\end{equation}
which is a product, over the row index $i\in[m]$, of ordinary multinomial coefficients.
If $A,B$ are fillings, and if we write $A\le B$, then we mean componentwise inequality, i.e., $A_{i,\alpha}\le B_{i,\alpha}$, for all $(i,\alpha)\in\mathbb{D}$.
For two fillings $A,B$ we will also denote by $AB$ their pointwise product, i.e.,
$[AB]_{i,\alpha}:=A_{i,\alpha}B_{i,\alpha}$, for all $(i,\alpha)\in\mathbb{D}$. This will primarily be used in conjunction with the next notion of indicator function or filling.
If $\mathbb{S}$ is a subset of the diagram $\mathbb{D}$, we define the effective filling $\mathbf{1}_{\mathbb{S}}$ by letting, for all $(i,\alpha)\in\mathbb{D}$,
\[
[\mathbf{1}_{\mathbb{S}}]_{i,\alpha}:=\bbone\{(i,\alpha)\in\mathbb{S}\}\ .
\]
Similarly to Iverson's bracket, we will use $\bbone\{\cdots\}$ for the indicator function of the enclosed condition, namely, equal to 1 if the property holds, and 0 otherwise.
 
Picking up the thread from the introduction, we call a transversal any subset $\mathbb{E}\subset\mathbb{D}$ such that,
$\forall i\in[m]$, $|\mathbb{E}\cap\mathbb{D}_i|=1$. Namely, we pick exactly one cell in each row of the diagram.
Such a transversal has the form
\[
\mathbb{E}=\{(1,\varepsilon_{1\ast}),\ldots,(m,\varepsilon_{m\ast})\}\ ,
\]
where the multiindices $\vepi$ can be conveniently stored as the rows of an $m\times n$ matrix $\varepsilon$.
If $A$ is a filling, the {\it support} of $A$ is
\[
{\rm supp}(A):=\{(i,\alpha)\in\mathbb{D}\ |\ A_{i,\alpha}\neq 0\}\ .
\]
Finally, the last important tool we need is the {\it push} operator to a given transversal.
If $\mathbb{E}$ is a transversal, and $A$ is a filling, then the push of $A$ to $\mathbb{E}$ is the unique filling $P_{\mathbb{E}}(A)$ such that ${\rm Row}(P_{\mathbb{E}}(A))={\rm Row}(A)$ and ${\rm supp}(P_{\mathbb{E}}(A))\subset\mathbb{E}$. In other words,
\[
[P_{\mathbb{E}}(A)]_{i,\alpha}:={\rm Row}_i(A)\ \bbone\{(i,\alpha)\in\mathbb{E}\}\ ,
\]
for all cells $(i,\alpha)\in\mathbb{D}$. Clearly, this is an idempotent or projection operator.
In subsequent computations, it will be useful to remember
that the operations $|\cdot|$, $\|\cdot\|$, ${\rm wt}(\cdot)$, ${\rm Row}(\cdot)$, $P_{\mathbb{E}}(\cdot)$, acting on weight vectors or on fillings,
are all $\mathbb{Z}$-linear.

\medskip
Without further ado, we put all the above notation into practice, with the precise computation of the quantities denoted by $\overline{\mathscr{R}}^{\mathbb{E}}$ and $\widehat{\mathscr{R}}^{\mathbb{E}}$, when performing the MH substitution mentioned in the introduction.
We fix the format $(m,n,d)$, and consider $\mathscr{R}\in\mathbb{Q}[\ua]$ which is multihomogeneous in the (indeterminate) coefficients $(a_{i,\alpha})_{\alpha\in\widetilde{\mathbb{D}}_i}$ of the forms $F_i$, $1\le i\le m$,
of multidegree $(\delta_1,\ldots,\delta_m)\in\mathbb{N}^m$, in these $m$ subcollections of $a$ variables. We will use the {\it column} vector $\delta=(\delta_1,\ldots,\delta_m)^{\rm T}$ to most conveniently store multidegree information.
We now have
\[
\mathscr{R}=\sum_{A\in {\rm EF}(\delta)} r_A\ \ua^{A}\ ,
\]
where
\[
{\rm EF}(\delta):=\{A\in\mathbb{N}^{\mathbb{D}}\ |\ {\rm Row}(A)=\delta\}\ ,
\]
and the coefficients $r_A$ are simply rational numbers.
After performing the MH substitution (\ref{MHsub}) relative to the given transversal $\mathbb{E}$ encoded in the $m\times n$ matrix $\varepsilon$ with rows $\vepi$, $1\le i\le m$, the polynomial
\[
\mathscr{R}=\sum_{A\in {\rm EF}(\delta)} r_A\ \ua^{\mathbf{1}_{\mathbb{D}\backslash\mathbb{E}}A}\times \prod_{i=1}^{m} a_{i,\vepi}^{A_{i,\vepi}}
\]
becomes
\[
\overline{\mathscr{R}}^{\mathbb{E}}:=
\sum_{A\in {\rm EF}(\delta)} r_A\ \ua^{\mathbf{1}_{\mathbb{D}\backslash\mathbb{E}}A}\times \prod_{i=1}^{m} 
\left(-\ \frac{\widetilde{F}_i}{\ux^{\vepi}}\right)^{A_{i,\vepi}}\ .
\]
The denominator is
\[
\prod_{i=1}^{m} 
\left(\ux^{\vepi}\right)^{A_{i,\vepi}}
=\ux^{\sum_{i=1}^{m}A_{i,\vepi} \vepi}
=\ux^{{\rm wt}(\mathbf{1}_{\mathbb{E}}A)}\ .
\]
Since $A$ is effective and, for all $i\in[m]$, $\sum_{\alpha\in\widetilde{\mathbb{D}}_i}A_{i,\alpha}=\delta_i$, we must have $0\le A_{i,\vepi}\le\delta_i$. Let ${\rm wt}_{\rm max}:=\sum_{i=1}^{m}\delta_i\vepi=\delta^{\rm T}\varepsilon$. Since the entries of $\vepi$ are nonnegative, we must have the componentwise inequality between weight vectors ${\rm wt}(\mathbf{1}_{\mathbb{E}}A)\le {\rm wt}_{\rm max}$, i.e., the worst denominator is $\ux^{{\rm wt}_{\rm max}}$. We thus go ahead and define
\[
\widehat{\mathscr{R}}^{\mathbb{E}}:=\ux^{{\rm wt}_{\rm max}}\times \overline{\mathscr{R}}^{\mathbb{E}}\ .
\]
As a result,
\begin{align*}
\widehat{\mathscr{R}}^{\mathbb{E}} & =  
\sum_{A\in {\rm EF}(\delta)} r_A\ \ua^{\mathbf{1}_{\mathbb{D}\backslash\mathbb{E}}A}
\ \ux^{\sum_{i=1}^{m}(\delta_i-A_{i,\vepi})\vepi}
\ (-1)^{\sum_{i=1}^{m}A_{i,\vepi}}
\ \prod_{i=1}^{m} 
\left(
\sum_{\alpha\in\widetilde{\mathbb{D}}_i\backslash\{\vepi\}}
a_{i,\alpha}\ \ux^{\alpha}
\right)^{A_{i,\vepi}} \\
 & = 
\sum_{A\in {\rm EF}(\delta)} r_A\ 
\ua^{\mathbf{1}_{\mathbb{D}\backslash\mathbb{E}}A}
\ \ux^{{\rm wt}(P_{\mathbb{E}}(\mathbf{1}_{\mathbb{D}\backslash\mathbb{E}}A))}
\ (-1)^{\|\mathbf{1}_{\mathbb{E}}A\|}
\ \prod_{i=1}^{m} 
\left(\sum_{\alpha\in\widetilde{\mathbb{D}}_i\backslash\{\vepi\}}
a_{i,\alpha}\ \ux^{\alpha}
\right)^{A_{i,\vepi}}
\ .
\end{align*}
Next, we use the multinomial theorem in order to expand each factor
\[
\left(\sum_{\alpha\in\widetilde{\mathbb{D}}_i\backslash\{\vepi\}}
a_{i,\alpha}\ \ux^{\alpha}\right)^{A_{i,\vepi}} \ ,
\]
and then take the expansion of the product of these last multinomial expansions, with the outcome
\begin{align*}
\prod_{i=1}^{m}\left(\sum_{\alpha\in\widetilde{\mathbb{D}}_i\backslash\{\vepi\}}
a_{i,\alpha}\ \ux^{\alpha}\right)^{A_{i,\vepi}} 
 & =  \ \ \sum_{M}
\begin{bmatrix}
{\rm Row}(\mathbf{1}_{\mathbb{E}}A) \\
M
\end{bmatrix}
\prod_{i=1}^{m}
\left(
\prod_{\alpha\in\widetilde{\mathbb{D}}_i\backslash\{\vepi\}}
(a_{i,\alpha}\ \ux^{\alpha})^{M_{i,\alpha}}
\right) \\
 & =  \ \ \sum_{M}
\begin{bmatrix}
{\rm Row}(\mathbf{1}_{\mathbb{E}}A) \\
M
\end{bmatrix}
\ua^M\ \ux^{{\rm wt}(M)}\ ,
\end{align*}
where the sum is over effective fillings $M$ such that ${\rm supp}(M)\subset\mathbb{D}\backslash\mathbb{E}$ and ${\rm Row}(M)={\rm Row}(\mathbf{1}_{\mathbb{E}}A)$.
Therefore,
\[
\widehat{\mathscr{R}}^{\mathbb{E}}=
\sum_{A,M} r_A
\ \ua^{\mathbf{1}_{\mathbb{D}\backslash\mathbb{E}}A+M}
\ (-1)^{\|\mathbf{1}_{\mathbb{E}}A\|}
\ \ux^{{\rm wt}(P_{\mathbb{E}}(\mathbf{1}_{\mathbb{D}\backslash\mathbb{E}}A))+{\rm wt}(M)}
\ \begin{bmatrix}
{\rm Row}(\mathbf{1}_{\mathbb{E}}A) \\
M
\end{bmatrix}\ ,
\]
where the sum is over $A\in{\rm EF}(\delta)$ and $M\in{\rm EF}({\rm Row}(\mathbf{1}_{\mathbb{E}}A))$ with ${\rm supp}(M)\subset\mathbb{D}\backslash\mathbb{E}$.
Given an effective filling $C\in\mathbb{N}^{\mathbb{D}}$, and a multiindex $\beta\in\mathbb{N}^n$, we can make precise the equation (\ref{shortmaineq})
denoted earlier by ${\rm Eq}(\mathbb{E},C,\beta)$:
\begin{equation}
\sum_{\substack{A\in{\rm EF}(\delta) \\
M\in{\rm EF}({\rm Row}(\mathbf{1}_{\mathbb{E}}A)) }} 
\bbone\left\{
\begin{array}{c}
\mathbf{1}_{\mathbb{D}\backslash\mathbb{E}}A+M=C \\
{\rm wt}(P_{\mathbb{E}}(\mathbf{1}_{\mathbb{D}\backslash\mathbb{E}}A))+{\rm wt}(M)=\beta \\
{\rm supp}(M)\subset\mathbb{D}\backslash\mathbb{E}
\end{array}
\right\}
(-1)^{\|\mathbf{1}_{\mathbb{E}}A\|}
\begin{bmatrix}
{\rm Row}(\mathbf{1}_{\mathbb{E}}A) \\
M
\end{bmatrix}
r_A
=0\ .
\label{longmaineq}
\end{equation}
While the range of the parameter $\mathbb{E}$ for the above equation is clear, namely, the set of all transversals for the given format, let us now examine the ranges of $C$ and $\beta$. Since the equation arises from extracting coefficients from a polynomial in $\mathbb{Q}[\ua_{\mathbb{D}\backslash\mathbb{E}},\ux]$, we must have ${\rm supp}(C)\subset\mathbb{D}\backslash\mathbb{E}$. Assuming the equation does not reduce to the vacuous statement $0=0$, i.e., there exist $A$ and $M$ satisfying the requirements in (\ref{longmaineq}), we can deduce the following constraints. We must have
\begin{eqnarray*}
{\rm Row}(C) & = & {\rm Row}\left(\mathbf{1}_{\mathbb{D}\backslash\mathbb{E}}A+M\right) \\
 & = & {\rm Row}\left(A-\mathbf{1}_{\mathbb{E}}A+M\right) \\
 & = & \delta -{\rm Row}\left(\mathbf{1}_{\mathbb{E}}A\right)
+{\rm Row}\left(M\right) \\
 & = & \delta\ .
\end{eqnarray*}
We must also have, using the weight condition on $A$ and $M$ in (\ref{longmaineq}), the identity (\ref{obviouswteq}), and the preservation of row sum vectors by push operators,
\begin{eqnarray*}
|\beta| & = & \left|{\rm wt}\left(P_{\mathbb{E}}(\mathbf{1}_{\mathbb{D}\backslash\mathbb{E}}A)+M\right) \right| \\
 & = & d^{\rm T}{\rm Row}\left(P_{\mathbb{E}}(\mathbf{1}_{\mathbb{D}\backslash\mathbb{E}}A)+M\right) \\
 & = & d^{\rm T}\left(
{\rm Row}\left(\mathbf{1}_{\mathbb{D}\backslash\mathbb{E}}A\right)+{\rm Row}(\mathbf{1}_{\mathbb{E}}A)
\right) \\
 & = & d^{\rm T}\delta\ .
\end{eqnarray*}
As a summary of the previous discussion, we have established the following proposition.

\begin{proposition}
\label{ifdzcharprop}
Let $(m,n,d)$ be a format and $\mathbb{E}$ be a transversal for that format. Let $\delta^{\rm T}=(\delta_1,\ldots,\delta_m)\in\mathbb{N}^m$. A polynomial $\mathscr{R}=
\sum_{A}r_A\ \ua^A\in\mathbb{Q}[\ua]$ which is multihomogeneous of multidegree $\delta^{\rm T}$, is such that $\widehat{\mathscr{R}}^{\mathbb{E}}$ is identically zero
iff, for all $C\in{\rm EF}(\delta)$ such that ${\rm supp}(C)\subset\mathbb{D}\backslash\mathbb{E}$, and for all $\beta\in\mathbb{N}^n$ such that $|\beta|=d^{\rm T}\delta$, we have that ${\rm Eq}(\mathbb{E},C,\beta)$ holds.
\end{proposition}

From now on, when referring to an equation ${\rm Eq}(\mathbb{E},C,\beta)$, we will assume $\mathbb{E},C,\beta$ satisfy the conditions stated in Proposition \ref{ifdzcharprop}. Note that the double sum over $A$ and $M$ is redundant since there is a one-to-one correspondence between these `summation indices'. We thus have two variants of ${\rm Eq}(\mathbb{E},C,\beta)$, one as a sum over $A$, and another one as a sum over $M$. 
Indeed, from the condition $\mathbf{1}_{\mathbb{D}\backslash\mathbb{E}}A+M=C $ in (\ref{longmaineq}), one trivially recovers $M$ from $A$ by
\[
M=C-\mathbf{1}_{\mathbb{D}\backslash\mathbb{E}}A\ . 
\]
Conversely, we note that ${\rm Row}(M)={\rm Row}(\mathbf{1}_{\mathbb{E}}A)$, and ${\rm supp}(\mathbf{1}_{\mathbb{E}}A)\subset\mathbb{E}$ imply $\mathbf{1}_{\mathbb{E}}A=P_{\mathbb{E}}(M)$, by definition of the push operator.
We can thus recover $A$ from $M$ by
\[
A=\mathbf{1}_{\mathbb{D}\backslash\mathbb{E}}A+\mathbf{1}_{\mathbb{E}}A
=C-M+P_{\mathbb{E}}(M)\ .
\]
Note that, because $P_{\mathbb{E}}\left(\mathbf{1}_{\mathbb{E}}A\right)=\mathbf{1}_{\mathbb{E}}A$, we have
\begin{eqnarray*}
{\rm wt}\left(P_{\mathbb{E}}(\mathbf{1}_{\mathbb{D}\backslash\mathbb{E}}A)\right)
& = & {\rm wt}\left(P_{\mathbb{E}}(A-\mathbf{1}_{\mathbb{E}}A)\right) \\
 & = &  {\rm wt}\left(P_{\mathbb{E}}(A)\right)-{\rm wt}\left(\mathbf{1}_{\mathbb{E}}A\right) \\
& = & \sum_{i=1}^{m}{\rm Row}_i(A)\vepi -{\rm wt}\left(\mathbf{1}_{\mathbb{E}}A\right) \\
& = & \delta^{\rm T}\varepsilon-{\rm wt}\left(\mathbf{1}_{\mathbb{E}}A\right) \ .
\end{eqnarray*}
Since ${\rm wt}(M)={\rm wt}(C)-{\rm wt}\left(\mathbf{1}_{\mathbb{D}\backslash\mathbb{E}}A\right)$, we see that the weight condition in (\ref{longmaineq}) becomes, after eliminating $M$, the condition
\[
{\rm wt}(A)={\rm wt}(C)-\beta+\delta^{\rm T}\varepsilon\ .
\]
We thus see that our equations are {\it isobaric}, namely, all the $A$'s featuring in such an equation must have the same weight vector ${\rm wt}(A)$.
We can now easily write the $A$ variant of ${\rm Eq}(\mathbb{E},C,\beta)$ which is an equivalent reformulation of (\ref{longmaineq}):
\begin{equation}
\sum_{A\in{\rm EF}(\delta)} 
\bbone\left\{
\begin{array}{c}
\mathbf{1}_{\mathbb{D}\backslash\mathbb{E}}A\le C \\
{\rm wt}(A)={\rm wt}(C)-\beta+\delta^{\rm T}\varepsilon
\end{array}
\right\}
(-1)^{\|\mathbf{1}_{\mathbb{E}}A\|}
\begin{bmatrix}
{\rm Row}(\mathbf{1}_{\mathbb{E}}A) \\
C-\mathbf{1}_{\mathbb{D}\backslash\mathbb{E}}A
\end{bmatrix}
r_A
=0\ .
\label{longmaineqA}
\end{equation}
From the remark $\mathbf{1}_{\mathbb{E}}A=P_{\mathbb{E}}(M)$, we can easily deduce
\[
{\rm wt}\left(P_{\mathbb{E}}(\mathbf{1}_{\mathbb{D}\backslash\mathbb{E}}A)\right)
=\delta^{\rm T}\varepsilon-{\rm wt}\left(P_{\mathbb{E}}(M)\right)\ ,
\]
and $\|\mathbf{1}_{\mathbb{E}}A\|=\|M\|$. We can then readily eliminate $A$ and write the $M$ variant of ${\rm Eq}(\mathbb{E},C,\beta)$ which is also an equivalent reformulation of (\ref{longmaineq}):
\begin{equation}
\sum_{M\in\mathbb{N}^{\mathbb{D}}} 
\bbone\left\{
\begin{array}{c}
M\le C \\
{\rm wt}\left(P_{\mathbb{E}}(M)\right)-{\rm wt}(M)=\delta^{\rm T}\varepsilon -\beta
\end{array}
\right\}
(-1)^{\|M\|}
\begin{bmatrix}
{\rm Row}(M) \\
M
\end{bmatrix}
r_{C-M+P_{\mathbb{E}}(M)}
=0\ .
\label{longmaineqM}
\end{equation}
Note that the condition ${\rm supp}(M)\subset\mathbb{D}\backslash\mathbb{E}$ is taken care of by the standing assumption ${\rm supp}(C)\subset\mathbb{D}\backslash\mathbb{E}$ and the requirement $M\le C$.
The $M$ variant will be better suited to understanding the mechanism behind why extremal pairs of transversals work, in \S\ref{warmupsec}, whereas the $A$ variant will be more convenient for the resolution of the system for the $r_A$ coefficients, in \S\ref{formulasec}.

In what follows we assume the format $(m,n,d)$ is fixed, and the polynomial $\mathscr{R}\in\mathbb{Q}[\ua]$ is multihomogeneous of multidegree $\delta^T$.
Given a multiindex $\gamma\in\mathbb{N}^n$, we say that $\mathscr{R}$ satisfies the property B\'ezout($\gamma$) if there exists polynomials $G_1,\ldots,G_m$ in $\mathbb{Q}[\ua,\ux]$ for which the relation (\ref{Bezout}) holds.

\begin{proposition}\label{FromBezoutProp}
If $\mathscr{R}$ satisfies  B\'ezout($\gamma$), for some  $\gamma\in\mathbb{N}^n$, then for all transversal $\mathbb{E}$, we have $\widehat{\mathscr{R}}^{\mathbb{E}}=0$.
\end{proposition}

\noindent{\bf Proof:}
By design, the MH substitution forces $F_i=a_{i,\vepi}\ \ux^{\vepi}+\widetilde{F}_i$ to become zero in the ring
$\mathbb{Q}[\ux,\ux^{-1}][\ua_{\mathbb{D}\backslash\mathbb{E}}]$, and the right-hand side of the B\'ezout identity (\ref{Bezout}) thus vanishes. This forces $\overline{\mathscr{R}}^{\mathbb{E}}$ on the left-hand side to also vanish, and after multiplication by $\ux^{{\rm wt}_{\rm max}}$, we also get $\widehat{\mathscr{R}}^{\mathbb{E}}=0$ in the ring $\mathbb{Q}[\ua,\ux]$.
\qed

\begin{proposition}\label{ToBezoutProp}
If the transversal $\mathbb{E}$ is such that $\widehat{\mathscr{R}}^{\mathbb{E}}=0$, then $\mathscr{R}$ satisfies B\'ezout($\gamma$) with $\gamma=\delta^{\rm T}\varepsilon$.
\end{proposition}

\noindent{\bf Proof:}
For the following computation, we will be working over the ring $\mathbb{Q}[\ux,\ux^{-1}][\ua_{\mathbb{D}}]$. Recall that, for all $i\in[m]$, $F_i:=\sum_{\alpha\in\widetilde{\mathbb{D}}_i}a_{i,\alpha}\ \ux^{\alpha}$, whereas
$\widetilde{F}_i:=\sum_{\alpha\in\widetilde{\mathbb{D}}_i\backslash\{\vepi\}}a_{i,\alpha}\ \ux^{\alpha}$. Assuming $\mathscr{R}$ is given by the expansion $\mathscr{R} = \sum_{A\in{\rm EF}(\delta)} r_A\ \ua^A $, we can write
\begin{eqnarray*}
\mathscr{R} & = &  \sum_{A\in{\rm EF}(\delta)} r_A\ \ua^{\mathbf{1}_{\mathbb{D}\backslash\mathbb{E}}A}\times\prod_{i=1}^{m}
a_{i,\vepi}^{A_{i,\vepi}} \\
 & = &  \sum_{A\in{\rm EF}(\delta)} r_A\ \ua^{\mathbf{1}_{\mathbb{D}\backslash\mathbb{E}}A}\times\prod_{i=1}^{m}
\left(
\frac{{F}_i}{\ux^{\vepi}}-\ \frac{\widetilde{F}_i}{\ux^{\vepi}}
\right)^{A_{i,\vepi}} \\
 & = &  \sum_{A\in{\rm EF}(\delta)} r_A\ \ua^{\mathbf{1}_{\mathbb{D}\backslash\mathbb{E}}A}\times\prod_{i=1}^{m}
\left(
\sum_{k_i=0}^{A_{i,\vepi}}\binom{A_{i,\vepi}}{k_i}
\left(\frac{{F}_i}{\ux^{\vepi}}\right)^{k_i}
\left(- \frac{\widetilde{F}_i}{\ux^{\vepi}}\right)^{A_{i,\vepi}-k_i}
\right) \\
 & = & 
\overline{\mathscr{R}}^{\mathbb{E}}
+\sum_{A\in{\rm EF}(\delta)}
\sum_{\substack{k\in\mathbb{N}^m\backslash\{0\} \\ \forall i,\ k_i\le A_{i,\vepi}}}
r_A\times 
\frac{\ua^{\mathbf{1}_{\mathbb{D}\backslash\mathbb{E}}A}}{\ux^{{\rm wt}(\mathbf{1}_{\mathbb{E}}A)}}
\times\prod_{i=1}^{m}
\left(
\binom{A_{i,\vepi}}{k_i}\ F_i^{k_i}\ (-\widetilde{F}_i)^{A_{i,\vepi}-k_i}
\right)\ ,
\end{eqnarray*}
where $\overline{\mathscr{R}}^{\mathbb{E}}$, previously defined in $\mathbb{Q}[\ux,\ux^{-1}][\ua_{\mathbb{D}\backslash\mathbb{E}}]$,
is now seen as an element of the larger ring $\mathbb{Q}[\ux,\ux^{-1}][\ua_{\mathbb{D}}]$.
We multiply the last equation by $\ux^{\delta^{\rm T}\varepsilon}$, and note that $\ux^{\delta^{\rm T}\varepsilon}
\ \overline{\mathscr{R}}^{\mathbb{E}}=\widehat{\mathscr{R}}^{\mathbb{E}}=0$,
by assumption.
We then immediately get the wanted B\'ezout relation (\ref{Bezout}) in $\mathbb{Q}[\ua,\ux]$, with, e.g., the choice of polynomials
\begin{eqnarray*}
\lefteqn{
G_q=\sum_{A\in{\rm EF}(\delta)}
\ \ \sum_{k\in\mathbb{N}^m\backslash\{0\}}
\bbone\left\{
\begin{array}{c}
\forall i, k_i\le A_{i,\vepi} \\
k_q\ge 1 \\
\forall i<q, k_i=0
\end{array}
\right\}\times\qquad\qquad} & & \\
 & & \qquad\qquad r_A\ \ua^{\mathbf{1}_{\mathbb{D}\backslash\mathbb{E}}A}
\ \ux^{\delta^{\rm T}\varepsilon-{\rm wt}(\mathbf{1}_{\mathbb{E}}A)}
\times\prod_{i=1}^{m}
\left(
\binom{A_{i,\vepi}}{k_i}\ F_i^{k_i-\bbone\{i=q\}}\ (-\widetilde{F}_i)^{A_{i,\vepi}-k_i}
\right)
\ ,
\end{eqnarray*}
for $1\le q\le m$.
\qed

\begin{proposition}\label{TFAEprop}
The following statements are equivalent.
\begin{enumerate}[label=(\roman*)]
\item
$\mathscr{R}$ is an IFDZ, i.e., there exists $q\in\mathbb{N}$, such that for all multiindex $\gamma\in\mathbb{N}^n$ with $|\gamma|=q$, property B\'ezout($\gamma$) holds for $\mathscr{R}$.
\item
There exists $\gamma\in\mathbb{N}^n$, such that property B\'ezout($\gamma$) holds for $\mathscr{R}$.
\item
For all transversal $\mathbb{E}$, we have $\widehat{\mathscr{R}}^{\mathbb{E}}=0$.
\item
There exists a transversal $\mathbb{E}$, such that $\widehat{\mathscr{R}}^{\mathbb{E}}=0$.
\end{enumerate}
\end{proposition}

\noindent{\bf Proof:}
Proposition \ref{FromBezoutProp} shows the implication {\it(ii)}$\Rightarrow${\it(iii)}.
The implications {\it(i)}$\Rightarrow${\it(ii)} and {\it(iii)}$\Rightarrow${\it(iv)} are trivial. 
Proposition \ref{ToBezoutProp} shows the implication {\it(iv)}$\Rightarrow${\it(ii)}.
All that remains is to show that any one of {\it(ii)}, {\it(iii)}, or {\it(iv)} implies {\it(i)}. Suppose {\it(iii)} holds, and let $j\in[n]$.
Consider the transversal
\[
\mathbb{E}=\{(i,d_i e_j)\ |\ 1\le i\le m\}\ ,
\]
with its associated matrix $\varepsilon$, and corresponding monomial selection $x_j^{d_1},\ldots,x_j^{d_m}$.
Note that
\[
\delta^{\rm T}\varepsilon=(\delta_1,\ldots,\delta_m)
\begin{pmatrix}
0 & \cdots & 0 & d_1 & 0 & \cdots & 0 \\
0 & \cdots & 0 & d_2 & 0 & \cdots & 0 \\
\vdots &  & \vdots & \vdots & \vdots &  & \vdots \\
0 & \cdots & 0 & d_m & 0 & \cdots & 0 
\end{pmatrix}
=(\delta^{\rm T}d)\ e_j\ .
\]
If $\delta=0$, then $\mathscr{R}$ is constant and {\it(iii)} implies this constant is zero, in which case {\it(i)} holds trivially. Hence, we can assume $\delta\neq 0$ and, therefore, $\delta^{\rm T}d\ge 1$.
We have, by Proposition \ref{ToBezoutProp}, the existence of polynomials $G_{ij}\in\mathbb{Q}[\ua,\ux]$, $1\le i\le m$,
such that
\[
x_j^{\delta^{\rm T}d}\ \mathscr{R}=\sum_{i=1}^{m} F_i\ G_{ij}\ .
\]
Take $q=n(\delta^{\rm T}d)-n+1$, and suppose the multiindex $\gamma=(\gamma_1,\ldots,\gamma_n)$ is such that $|\gamma|=q$. Clearly, there exists $j\in[n]$ such that $\gamma_j\ge \delta^{\rm T}d$, and this allows us to write
\[
\ux^{\gamma}\ \mathscr{R}=\sum_{i=1}^{m} F_i \left(G_{ij}\ \ux^{\gamma-(\delta^{\rm T}d)e_j}\right) \ .
\]
Hence, B\'ezout($\gamma$) holds.
\qed

\begin{remark}\label{singlexErem}
The previous propositions are due to Mertens and Hurwitz (see the pedagogical presentations in~\cite[Ch. XI]{VdWaerden} and~\cite[Ch. IV]{HodgeP}), except perhaps for the scope of generality for the transversals $\mathbb{E}$. In these references, the transversals used are as in the proof of Proposition \ref{TFAEprop}, namely, the selected monomials are all powers of the same variable $x_j$. Apart from the linear case, we were not able to use such restricted transversals in order to solve the corresponding system of equations ${\rm Eq}(\mathbb{E},C,\beta)$. Allowing more general transversals like the one in (\ref{diagonaleps}) proved to be a more fruitful choice, as far as solving such systems.
\end{remark}

\begin{lemma}\label{extensionlem}
Let $m,q\in\mathbb{N}$ be such that $m\le q$, and consider a collection of integers $d_1,\ldots,d_q\ge 1$ and another one $\delta_1,\ldots,\delta_q\ge 0$.
Suppose there exists a nonzero IFDZ $\mathscr{R}$ for the format $(m,n,(d_1,\ldots,d_m)^{\rm T})$ and multidegree $(\delta_1,\ldots,\delta_m)$.
Then, there exists a nonzero IFDZ $\mathscr{S}$ for the format $(q,n,(d_1,\ldots,d_q)^{\rm T})$ and multidegree $(\delta_1,\ldots,\delta_q)$.
\end{lemma}

\noindent{\bf Proof:}
We see the diagram $\mathbb{D}$ corresponding to the format  $(m,n,(d_1,\ldots,d_m)^{\rm T})$ as the protion made of the first $m$ rows of the larger or extended diagram $\mathbb{D}^{\rm ext}$ corresponding to the format
$(q,n,(d_1,\ldots,d_q)^{\rm T})$. By hypothesis, and by Proposition \ref{TFAEprop}, we have a nonzero IFDZ $\mathscr{R}\in\mathbb{Q}[\ua_{\mathbb{D}}]$, a multiindex $\gamma\in\mathbb{N}^n$, and a collection of polynomials $G_1,\ldots,G_m\in\mathbb{Q}[\ua_{\mathbb{D}},\ux]$
such that
\[
\ux^{\gamma}\ \mathscr{R}=F_1\ G_1+\cdots+ F_m\ G_m\ ,
\]
in the ring $\mathbb{Q}[\ua_{\mathbb{D}},\ux]$.
We add new homogeneous polynomials $F_i$, $m< i\le q$
with their indeterminate coefficients $a_{i,\alpha}$, so the original collection of formal variables $\ua_{\mathbb{D}}=(a_{i,\alpha})_{(i,\alpha)\in\mathbb{D}}$
is seen as a subcollection of $\ua_{\mathbb{D}^{\rm ext}}=(a_{i,\alpha})_{(i,\alpha)\in\mathbb{D}^{\rm ext}}$.
We will also use the ring inclusions $\mathbb{Q}[\ua_{\mathbb{D}}]\subset\mathbb{Q}[\ua_{\mathbb{D}^{\rm ext}}]$ and $\mathbb{Q}[\ua_{\mathbb{D}},\ux]\subset\mathbb{Q}[\ua_{\mathbb{D}^{\rm ext}},\ux]$.
Let
\[
\mathscr{S}=\mathscr{R}\times \mathscr{A}\ ,
\]
where
\[
\mathscr{A}=(a_{m+1,\ d_{m+1}e_1})^{\delta_{m+1}}
\cdots (a_{q,\ d_q e_1})^{\delta_q}\ .
\]
By multiplying the previous B\'ezout relation, seen in the ring $\mathbb{Q}[\ua_{\mathbb{D}^{\rm ext}},\ux]$, by $\mathscr{A}$, we get
\[
\ux^{\gamma}\mathscr{S}=F_1 (G_1\mathscr{A})+\cdots+ F_m (G_m\mathscr{A})
+F_{m+1}\times 0+\cdots +F_q\times 0\ .
\]
Clearly, $\mathscr{S}\neq 0$, and, by Proposition \ref{TFAEprop}, it is an IFDZ of format $(q,n,(d_1,\ldots,d_q)^{\rm T})$ and multidegree $(\delta_1,\ldots,\delta_q)$.
\qed

\begin{theorem}\label{lowdegthm}
Let $\mathscr{R}$ be an IFDZ for the format $(n,n,(d_1,\ldots,d_n)^{\rm T})$ and with multidegree $(\delta_1,\ldots,\delta_n)$.
If $\exists i\in[n], \delta_i<\prod_{j\in[n]\backslash\{i\}}d_j$, then $\mathscr{R}=0$.
\end{theorem}

The proof of this theorem is postponed to \S\ref{warmupsec}.

\begin{theorem}
Suppose $m<n$.
If $\mathscr{R}$ is an IFDZ for the format $(m,n,(d_1,\ldots,d_m)^{\rm T})$ and with multidegree $(\delta_1,\ldots,\delta_m)$, then $\mathscr{R}=0$.
\end{theorem}

\noindent{\bf Proof:}
Arguing by contradiction, suppose $\mathscr{R}\neq 0$. We let $q=n$, $d_{m+1}=\cdots=d_n=1$, and $\delta_{m+1}=\cdots=\delta_n=0$, and use Lemma \ref{extensionlem} in order to produce a nonzero IFDZ for the format $(n,n,(d_1,\ldots,d_n)^{\rm T})$ and with multidegree $(\delta_1,\ldots,\delta_n)$.
By Theorem \ref{lowdegthm}, we must have $0=\delta_n\ge d_1\cdots d_{n-1}=d_1\cdots d_m\ge 1$ which is a contradiction. Thus $\mathscr{R}$ must vanish.
\qed

\begin{remark}
The last theorem is not new (see~\cite[Prop. 4.3]{Jouanolou}) and is attributed to Hurwitz who proved it by a clever induction on the dimension, which includes the acyclicity of the Koszul complex for generic forms, in degree 1, as part of the induction hypothesis. See also~\cite[p. 7]{MacaulayFormulas}, regarding this acyclicity result.
\end{remark}

\section{Warm-up practice with the main reduction algorithm}
\label{warmupsec}

In the first part of this section we will prove Theorem \ref{lowdegthm}, as an introductory example of application for our new method. We assume $\mathscr{R}=\sum_{A\in {\rm EF}(\delta)} r_A\ \ua^{A}$ is an IFDZ for the format $(n,n,d)$, with $d=(d_1,\ldots,d_n)^{\rm T}$, and with multidegree $\delta^{\rm T}$. We assume 
$\exists i\in[n], \delta_i<\prod_{j\in[n]\backslash\{i\}}d_j$. Up to reordering the forms $F_1,\ldots,F_n$, there is no loss in generality in only considering the case $i=1$, i.e., assuming
\begin{equation}
\delta_1<d_2\cdots d_n\ .
\label{delta1assum}
\end{equation}
We pick the transversal $\mathbb{E}$ given by the matrix $\varepsilon$ in (\ref{diagonaleps}). This transversal $\mathbb{E}$ will be fixed throughout the following proof.
Since $\mathscr{R}$ is an IFDZ, by Proposition \ref{TFAEprop}, part {\it(iii)}, we have that $\widehat{\mathscr{R}}^{\mathbb{E}}=0$, for our chosen transversal. By Proposition \ref{ifdzcharprop}, this means the rational numbers $r_A$, $A\in {\rm EF}(\delta)$, must satisfy all the equations ${\rm Eq}(\mathbb{E},C,\beta)$, where $C\in{\rm EF}(\delta)$, with ${\rm supp}(C)\subset\mathbb{D}\backslash\mathbb{E}$, and $\beta\in\mathbb{N}^n$, with $|\beta|=d^{\rm T}\delta$.

Let $A$ be some effective filling in ${\rm EF}(\delta)$. We now explain our prescription for choosing the resolvent equation ${\rm Eq}(\mathbb{E},C,\beta)$ for the unknown $r_A$. Recall, from the introduction, that we are looking for equations where the coefficient of $r_A$ is $\pm 1$. In the $M$ representation (\ref{longmaineqM}), this means that the associated filling $M$ must be such that
\[
\begin{bmatrix}
{\rm Row}(M) \\
M
\end{bmatrix}=1\ .
\]
Since this is a product, over the rows of the diagram $\mathbb{D}$, of ordinary multinomial coefficients, the latter must then be of the form
\[
\binom{q}{0,\ldots,0,q,0,\ldots,0}\ .
\]
In other words, there must exist a transversal $\mathbb{T}$, such that ${\rm supp}(M)\subset\mathbb{T}$. We go ahead and pick $\mathbb{T}$ as the transversal given by the associated matrix $\tau$ in (\ref{taupick}). This transversal will also be fixed throughout the rest of the proof. Note that, as subsets of $\mathbb{D}$, the transversals $\mathbb{E}$ and $\mathbb{T}$ are disjoint.
We now decide to successively choose:
\begin{eqnarray}
M & := & P_{\mathbb{T}}\left(\mathbf{1}_{\mathbb{E}}A\right) \nonumber \\
C & := & \mathbf{1}_{\mathbb{D}\backslash\mathbb{E}}A+M \nonumber \\
\beta & := & 
{\rm wt}(P_{\mathbb{E}}(\mathbf{1}_{\mathbb{D}\backslash\mathbb{E}}A))+{\rm wt}(M)\ .
\label{resolvpick}
\end{eqnarray}
The above is {\it how we choose the resolvent equation} ${\rm Eq}(\mathbb{E},C,\beta)$ for the unknown $r_A$. It is easy to see that $C$ and $\beta$ are in the allowed ranges. Note that a more precise choice of notation, in the following, would be to use $M_A$, $C_A$, $\beta_A$, instead of $M$, $C$, $\beta$, in order to emphasize that these quantities are fixed and determined by $A$. However, we prefer to err on the side of lighter notation.
By construction, this resolvent equation, in the variant (\ref{longmaineqM}), reads
\[
(-1)^{\|\mathbf{1}_{\mathbb{E}}A\|} r_A+
\sum_{\wtM\in\mathbb{N}^{\mathbb{D}}} 
\bbone\left\{
\begin{array}{c}
\wtM\neq M \ ,\ \wtM\le C \\
{\rm wt}\left(P_{\mathbb{E}}(\wtM)\right)-{\rm wt}(\wtM)=
{\rm wt}\left(P_{\mathbb{E}}(M)\right)-{\rm wt}(M)
\end{array}
\right\}\times
\qquad\qquad\qquad\qquad\qquad
\]
\begin{equation}
\qquad\qquad\qquad\qquad\qquad\qquad\qquad\qquad
(-1)^{\|\wtM\|}
\begin{bmatrix}
{\rm Row}(\wtM) \\
\wtM
\end{bmatrix}
r_{C-\wtM+P_{\mathbb{E}}(\wtM)}
=0\ ,
\label{resolventeq}
\end{equation}
where we had to rename the summation index $\wtM$ in order to distinguish it from $M$ which is determined by $A$.
For any {\it effective} filling $B$, we define the {\it counter}
\[
\mathcal{Z}(B):=\sum_{(i,\alpha)\in\mathbb{D}\backslash(\mathbb{E}\cup\mathbb{T})}B_{i,\alpha}=\|\mathbf{1}_{\mathbb{D}\backslash(\mathbb{E}\cup\mathbb{T})}B\|\ge 0\ .
\]
We will say that an effective filling $\wtA$ {\it features in the resolvent equation for} $A$ if $r_{\wtA}$ appears in (\ref{resolventeq}), namely, if one can write
$\wtA=C-\wtM+P_{\mathbb{E}}(\wtM)$ for some $\wtM\in\mathbb{N}^{\mathbb{D}}$ such that the {\it weight condition}
\begin{equation}
{\rm wt}\left(P_{\mathbb{E}}(\wtM)\right)-{\rm wt}(\wtM)=
{\rm wt}\left(P_{\mathbb{E}}(M)\right)-{\rm wt}(M)
\label{weighteq}
\end{equation}
holds, as well as the cell-wise inequality
\[
\wtM\le C\ .
\]
The following lemma is central to this entire article.

\begin{lemma}\label{keylem}
If $\wtA\neq A$ features in the resolvent equation for $A$, then $\mathcal{Z}(\wtA)<\mathcal{Z}(A)$.
\end{lemma}

\noindent{\bf Proof:} Let $\wtA$ be as in the statement of the lemma, and let $\wtM$ be the associated filling by the one-to-one correspondence discussed right after Proposition \ref{ifdzcharprop}. By the specification of the resolvent equation indicated in (\ref{resolvpick}), we have that
\begin{eqnarray}
\wtA & = & C-\wtM+P_{\mathbb{E}}(\wtM) \nonumber \\
 & = & \mathbf{1}_{\mathbb{D}\backslash\mathbb{E}}A+M-\wtM+P_{\mathbb{E}}(\wtM) \nonumber \\
 & = & \mathbf{1}_{\mathbb{T}}A+\mathbf{1}_{\mathbb{D}\backslash(\mathbb{E}\cup\mathbb{T})}A+P_{\mathbb{T}}\left(\mathbf{1}_{\mathbb{E}}A\right)
-\wtM+P_{\mathbb{E}}(\wtM)\ .
\label{wtAdecomp}
\end{eqnarray}
Keeping in mind the support restrictions on the terms of the previous equation, and the fact that $\mathbb{E}$, $\mathbb{T}$, $\mathbb{D}\backslash(\mathbb{E}\cup\mathbb{T})$ form a set partition of the diagram $\mathbb{D}$,
we obtain
\[
\mathbf{1}_{\mathbb{D}\backslash(\mathbb{E}\cup\mathbb{T})}\wtA=
0+\mathbf{1}_{\mathbb{D}\backslash(\mathbb{E}\cup\mathbb{T})}A+0
-\mathbf{1}_{\mathbb{D}\backslash(\mathbb{E}\cup\mathbb{T})}\wtM+0\ ,
\]
after taking the cell-wise product of (\ref{wtAdecomp}) with the indicator filling $\mathbf{1}_{\mathbb{D}\backslash(\mathbb{E}\cup\mathbb{T})}$.
As a result,
\[
\mathcal{Z}(\wtA)=\mathcal{Z}(A)-\sum_{(i,\alpha)\in\mathbb{D}\backslash(\mathbb{E}\cup\mathbb{T})} \wtM_{i,\alpha}\ .
\]
Since $\wtM$ is effective, i.e., the entries $\wtM_{i,\alpha}$ are nonnegative, we see that the counter always decreases, at least weakly:
\[
\mathcal{Z}(\wtA)\le\mathcal{Z}(A)\ .
\]
We now take a closer look at the weight condition (\ref{weighteq}).
Since ${\rm supp}(\wtM)\subset\mathbb{D}\backslash\mathbb{E}$, we have
\[
{\rm wt}(\wtM)=\sum_{(i,\alpha)\in\mathbb{D}\backslash\mathbb{E}} \wtM_{i,\alpha}\ \alpha\ ,
\]
as well as
\[
{\rm wt}\left(P_{\mathbb{E}}(\wtM)\right)=\sum_{(i,\alpha)\in\mathbb{D}\backslash\mathbb{E}} \wtM_{i,\alpha}\ \vepi\ ,
\]
and therefore,
\[
{\rm wt}\left(P_{\mathbb{E}}(\wtM)\right)-{\rm wt}(\wtM)=
\sum_{(i,\alpha)\in\mathbb{D}\backslash\mathbb{E}} \wtM_{i,\alpha}\ (\vepi-\alpha)\ .
\]
In the particular case $\wtM=M=P_{\mathbb{T}}\left(\mathbf{1}_{\mathbb{E}}A\right)$, we have
\[
{\rm wt}\left(P_{\mathbb{E}}(M)\right)-{\rm wt}(M)=
\sum_{i=1}^{n}A_{i,\vepi}\ (\vepi-\taui)\ ,
\]
and the weight condition (\ref{weighteq}) can thus be written as
\[
\sum_{i=1}^{n}\wtM_{i,\taui}\ (\vepi-\taui)+
\sum_{(i,\alpha)\in\mathbb{D}\backslash(\mathbb{E}\cup\mathbb{T})} \wtM_{i,\alpha}\ (\vepi-\alpha)
=\sum_{i=1}^{n}A_{i,\vepi}\ (\vepi-\taui)\ ,
\]
where we separated the contributions of the cells in $\mathbb{T}$ from the rest, in the left-hand side.
We now introduce the row vector $y=(y_1,\ldots,y_n)$ defined by
\[
y_i:=\wtM_{i,\taui}-A_{i,\vepi}\ ,
\]
for all $i\in[n]$, as well as the matrix $H:=\varepsilon-\tau$. This gives us the reformulation of the weight condition (\ref{weighteq}) that will be most useful to our considerations:
\begin{equation}
yH+\sum_{(i,\alpha)\in\mathbb{D}\backslash(\mathbb{E}\cup\mathbb{T})} \wtM_{i,\alpha}\ (\vepi-\alpha)=0\ .
\label{usefulWC}
\end{equation}
Note that, for all $i\in[n]$,
$0\le A_{i,\vepi}\le \delta_i$, and $0\le A_{i,\taui}\le \delta_i$,
because $A\in{\rm EF}(\delta)$.
From
\[
0\le \wtM\le C=\mathbf{1}_{\mathbb{D}\backslash\mathbb{E}}A+P_{\mathbb{T}}(\mathbf{1}_{\mathbb{E}}A)\ ,
\]
and specializing to the cells of the transversal $\mathbb{T}$, we infer that, for all $i\in[n]$,
\[
0\le \wtM_{i,\taui}\le A_{i,\taui}+A_{i,\vepi}\ .
\]
As a result, we have the crucial inequalities
\begin{equation}
-\delta_i\le -A_{i,\vepi}\le y_i\le A_{i,\taui}\le \delta_i\ ,
\label{crucialineq}
\end{equation}
for all $i\in[n]$.

We now go back to the consideration of the effective filling $\wtA\neq A$ which features in the resolvent equation for $A$, and assume for the sake of contradiction that the counter does not drop, i.e., $\mathcal{Z}(\wtA)=\mathcal{Z}(A)$.
This implies that $\wtM_{i,\alpha}=0$ for all cells $(i,\alpha)$ in $\mathbb{D}\backslash(\mathbb{E}\cup\mathbb{T})$, and the weight condition reduces to
\[
yH=0\ .
\]
Recall that with our ansatz for the transversal pair $(\mathbb{E},\mathbb{T})$, the matrix $H$ is explicitly given by
\begin{equation}
H=
\begin{pmatrix}
1 & -1 & 0 & \cdots & 0 \\
-d_2+1 & d_2 & -1 & \ddots & \vdots \\
\vdots & 0 & \ddots & \ddots & 0 \\
-d_{n-1}+1 & \vdots & \ddots & d_{n-1} & -1 \\
-d_n & 0 & \cdots & 0 & d_n
\end{pmatrix}\ ,
\label{goodH}
\end{equation}
and its row sums vanish. Let $\widetilde{H}$ be the $n\times (n-1)$ matrix obtained from $H$ by removing the first column. The linear system $yH=0$ is equivalent to $y\widetilde{H}=0$ or, more explicitly,
\[
\left\{
\begin{array}{ccl}
y_1 & = & d_2 y_2 \\
y_2 & = & d_3 y_3 \\
 & \vdots & \\
y_{n-1} & = & d_n y_n\ ,
\end{array}
\right.
\] 
which is readily solved by
\[
y=y_n\left(d_2\cdots d_n,d_3\cdots d_n,\ldots, d_n,1\right)\ .
\]
By the $i=1$ case of (\ref{crucialineq}), and the hypothesis (\ref{delta1assum}), we must have
\[
|y_1|\le \delta_1<d_2\cdots d_n\ .
\]
Since $y_1=d_2\cdots d_n y_n$, we deduce $|y_n|<1$ and thus $y_n=0$, because these are integers. Hence, the entire vector $y$ vanishes, and this forces $\wtM=M$, i.e., $\wtA=A$ which is a contradiction. In conclusion, the counter must drop, i.e., $\mathcal{Z}(\wtA)<\mathcal{Z}(A)$.
\qed

We now finish the proof of Theorem \ref{lowdegthm} and establish that for all $A\in{\rm EF}(\delta)$, $r_A=0$, by complete induction on the counter $\mathcal{Z}(A)\ge 0$.
If $\mathcal{Z}(A)=0$, then Lemma \ref{keylem} implies that there is no $\wtA$ featuring in the resolvent equation for $A$, other than $A$ itself. Namely the equation reduces to $\pm r_A=0$, which establishes the base case of the induction.
For the induction step notice that, again by Lemma \ref{keylem}, the resolvent equation has the form (\ref{lambdarhoeq})
mentioned in the introduction, with $\lambda=\pm 1$, and where $\wtA\prec A$ is the strict partial order defined by the condition $\mathcal{Z}(\wtA)<\mathcal{Z}(A)$.
Equivalently, the associated weak partial order $\wtA\preceq A$ is defined by $\wtA=A$ or $\mathcal{Z}(\wtA)<\mathcal{Z}(A)$. Since, by the induction hypothesis, all the $r_{\wtA}$ vanish, so does $r_A$.
This concludes the proof of Theorem \ref{lowdegthm}.

\medskip
Let $\beta$ be a weight vector in $\mathbb{Z}^n$.
We will say that $\beta$ is {\it balanced} if it is an integer multiple of $(1,\ldots,1)$.

\begin{theorem}\label{balancedthm}
Let $\mathscr{R}$ be an IFDZ for the format $(n,n,(d_1,\ldots,d_n)^{\rm T})$ and with canonical multidegree $(\delta_1,\ldots,\delta_n)$ given by $\delta_i:=\prod_{j\in[n]\backslash\{i\}}d_j$, for all $i\in[n]$.
If $A\in{\rm EF}(\delta)$ is such that ${\rm wt}(A)$ is not balanced, then $r_A=0$.
\end{theorem}

\noindent{\bf Proof:} Consider a multiindex $w=(w_1,\ldots,w_n)\in\mathbb{N}^n$. If $A\in{\rm EF}(\delta)$ is such that ${\rm wt}(A)=w$, then, by (\ref{obviouswteq}), $|w|=\delta^{\rm T}d=nd_1\cdots d_n$. Saying $w$ is not balanced means that $w\neq d_1\cdots d_n(1,\ldots,1)$. This implies $\exists i\in[n]$, $w_i<d_1\cdots d_n$.
We need to show that if this is the case, and ${\rm wt}(A)=w$, then $r_A=0$. 
We first note, that up to reordering the variables $x_1,\ldots,x_n$, there is no loss of generality in assuming $i=1$, i.e., $w_1<d_1\cdots d_n$.
We will proceed as in the proofs of Theorem \ref{lowdegthm} and Lemma \ref{keylem}, with some minor modifications. Since the equations ${\rm Eq}(\mathbb{E},C,\beta)$ are isobaric, all the fillings $A$, $\wtA$ in the following discussion will be assumed to have weight $w$.
We use the same choice of transversals and resolvent equation as before.

We first prove that the statement in Lemma \ref{keylem} also holds in the present situation. Assume for the sake of contradiction that there is no drop in the counter, i.e., $\mathcal{Z}(\wtA)=\mathcal{Z}(A)$, then we again have 
$\wtM_{i,\alpha}=0$, for all $(i,\alpha)\in\mathbb{D}\backslash(\mathbb{E}\cup\mathbb{T})$. The weight condition (\ref{usefulWC}) again reduces to the system $y\widetilde{H}=0$, with the same notations as before. The inequalities (\ref{crucialineq}) now give
\[
|y_1|\le \delta_1=d_2\cdots d_n\ ,
\]
while solving the system $y\widetilde{H}=0$ gives $y_1=d_2\cdots d_n y_n$. Hence, $|y_n|\le 1$. Again, $y_n=0$ is excluded by $\wtA\neq A$, i.e., $\wtM\neq M$. So the integer $y_n$ can only take the values $y_n=-1$ or $y_n=1$.

\smallskip
\noindent{\bf 1st case:} $y_n=-1$.

\smallskip
This implies $y_1=-\delta_1$.
From (\ref{crucialineq}), we see that $-\delta_1\le -A_{1,\varepsilon_{1\ast}}\le -\delta_1$, namely, $A_{1,\varepsilon_{1\ast}}=\delta_1$.
As a result,
\[
w={\rm wt}(A)=\sum_{(i,\alpha)\in\mathbb{D}}A_{i,\alpha} \ \alpha\ge
A_{1,\varepsilon_{1\ast}}\ \varepsilon_{1\ast}=\delta_1 d_1 e_1\ .
\]
where the inequality of weight vectors is understood componentwise.
The above holds because the entries of $A$ and the components of the $\alpha$'s are nonnegative.
Looking at the first component, we obtain $w_1\ge d_1\cdots d_n$ which is a contradiction.

\smallskip

\noindent{\bf 2nd case:} $y_n=1$.

\smallskip
This implies $y_1=\delta_1$.
From (\ref{crucialineq}), we see that $\delta_1\le A_{1,\tau_{1\ast}}\le \delta_1$, namely, $A_{1,\tau_{1\ast}}=\delta_1$.
Since ${\rm Row}_1(A)=\delta_1$, this forces $A_{1,\varepsilon_{1\ast}}=0$.
As a result, $\wtM_{1,\tau_{1\ast}}=y_1+A_{1,\varepsilon_{1\ast}}=\delta_1$.
Using the relation (\ref{wtAdecomp}), evaluated on the cell $(1,\varepsilon_{1\ast})$, we get
\[
\wtA_{1,\varepsilon_{1\ast}}=\wtM_{1,\tau_{1\ast}}=\delta_1\ ,
\]
where the only nonzero contribution comes from the last term $P_{\mathbb{E}}(\wtM)$.
By the same reasoning as in the 1st case, now applied to $\wtA$ instead of $A$, we see, by looking at the first component of the weight of $\wtA$, that  $w_1\ge d_1\cdots d_n$ which is a contradiction.

\smallskip
In all cases, we arrive at a contradiction, and the counter must therefore drop.
Having established the analogue of Lemma \ref{keylem}, Theorem \ref{balancedthm} follows by the same induction as in the previous proof of Theorem \ref{lowdegthm}. The only difference is that, instead of running the induction on all the fillings in ${\rm EF}(\delta)$,  we only consider those with weight given by $w$.
\qed

\section{A convex analysis interlude}
\label{convexsec}

In this section, the format of interest will be $(n,n,d)$, with $n\ge 2$, and $\mathbb{E}$, $\mathbb{T}$ will denote two arbitrary disjoint transversals in the diagram $\mathbb{D}$. Their associated matrices are $\varepsilon$, $\tau$, respectively. Here, we are not only discussing the special ones defined in (\ref{diagonaleps}) and (\ref{taupick}).
We again let $H=\varepsilon-\tau$, which satisfies $H(1,\ldots,1)^{\rm T}=0$.
We will denote by $\Delta_{ij}(H)$ the minor determinant of the matrix obtained by removing the $i$-th row, and the $j$-th column of $H$.
We will say that $(\mathbb{E},\mathbb{T})$ is an extremal pair of transversals if
$\exists i,j,k\in[n]$, with $\Delta_{jk}(H)\neq 0$, and
\[
\left|\frac{\Delta_{ik}(H)}{\Delta_{jk}(H)}\right|=\prod_{\ell\in[n]\backslash\{i\}}d_{\ell}\ .
\]
By taking $i=1$, $j=n$, and $k=1$, and in view of (\ref{goodH}), we see that the transversal pair defined by (\ref{diagonaleps}) and (\ref{taupick}) is indeed extremal.
If $y=(y_1,\ldots,y_n)$ is a row vector, the linear system $yH=0$ is equivalent to $y\widetilde{H}=0$, where $\widetilde{H}$ is the matrix obtained
from $H$ be removing the $k$-th column. 
Provided $\Delta_{jk}(H)\neq 0$, we can solve for all the components of $y$ in terms of $y_j$. A simple use of Cramer's rule gives
\[
y_i=(-1)^{i+j}\ y_j\times \frac{\Delta_{ik}(H)}{\Delta_{jk}(H)}\ .
\]
The main idea in Lemma \ref{keylem} is to arrange for the amplifying factor relating $y_i$ to $y_j$ to be as large in absolute value as possible. Since $H$ has integer entries, the nonvanishing denominator implies $|\Delta_{jk}(H)|\ge 1$. As a consequence of Theorem \ref{substochthm} below, one always has 
\begin{equation}
|\Delta_{ik}(H)|\le\prod_{\ell\in[n]\backslash\{i\}}d_{\ell}\ ,
\label{minorineq}
\end{equation}
as explained at the end of this section.
Therefore, for extremal pairs, the amplifying factor is as large as possible.

\medskip
Let $s\in\mathbb{Z}_{>0}$, and consider a matrix $P=(P_{ij})_{1\le i,j\le s}$, with nonnegative real entries. It is called a substochastic matrix if, for all $i\in[s]$, $\sum_{j=1}^{s}P_{ij}\le 1$.

\begin{theorem}\label{substochthm}
Let $s\in\mathbb{Z}_{>0}$ and let $P$, $Q$ be two $s\times s$ substochastic matrices. Then, we have the determinantal inequality
$|{\rm det}(P-Q)|\le 1$.
\end{theorem}

\noindent{\bf Proof:} 
For given row vectors $v_1,\ldots,v_s$ in $\mathbb{R}^s$, we denote by
${\rm det}(v_1,\ldots,v_s)$ the determinant of the $s\times s$ matrix whose $i$-th row is given by $v_i$, for all $i\in[s]$. The statement to be established can be rephrased as showing that if all the $v_i$'s belong to $K_s:=\Delta_s+(-\Delta_s)$, then that determinant is no greater than 1 in absolute value. Here $\Delta_s$
denotes the simplex
\[
\Delta_s:=\{(x_1,\ldots,x_s)\in\mathbb{R}^s\ |\ \forall i\in[s], x_i\ge 0,\ {\rm and}\ \sum_{i=1}^{s}x_i\le 1\}\ .
\]
For a subset $X$ of $\mathbb{R}^s$, we write $-X:=\{-x\ |\ x\in X\}$. The sum of two subsets $X$, $Y$ of $\mathbb{R}^s$, is defined, as usual, by $X+Y:=\{x+y\ |\ x\in X, y\in Y\}$.
Denoting the convex hull of a set $X$ by ${\rm Conv}(X)$, we have that
\[
\Delta_s={\rm Conv}(\{0,e_1,\ldots,e_s\})\ ,
\]
where, for the purpose of this proof, the $e_i$ are the standard basis vectors of $\mathbb{R}^s$, instead of $\mathbb{R}^n$, as in the rest of the article.
Recall that for two finite subsets $X$, $Y$ of $\mathbb{R}^s$, the Minkowski sum of their convex hulls satisfies
\[
{\rm Conv}(X)+{\rm Conv}(Y)={\rm Conv}(X+Y)\ .
\]
Indeed, the $\supset$ inclusion is trivial, whereas the $\subset$ inclusion follows from the existence of nonnegative matrices (or contingency tables) with imposed compatible row and column marginals (see, e.g.,~\cite[Thm. 2.1.1]{Brualdi}).
We therefore easily obtain that $K_s={\rm Conv}(E_s)$, where $E_s$ is the finite set of size $s(s+1)$ made of the vectors $\pm e_i$, for $i\in[s]$, and the vectors $e_i-e_j$, for $i,j\in[s]$ with $i\neq j$. While multilinear optimization is in general difficult, here the underlying set has a cartesian product structure which allows us to reduce the proof of the inequality $|{\rm det}(v_1,\ldots,v_s)|\le 1$, for arbitrary vectors $v_i$ in $K_s$, to the case where these vectors are all in $E_s$. Indeed, let us enumerate the elements of $E_s$ as $w_1,\ldots,w_q$ with $q=s(s+1)$, and let us consider vectors $v_i$ that are in $K_s$. Then, there exists an $s\times q$ matrix $\lambda$ with nonnegative entries such that $\forall i\in[s]$, $\sum_{j=1}^{q}\lambda_{ij}=1$, and $v_i=\sum_{j=1}^{q}\lambda_{ij} w_j$. Expanding by multilinearity, we have
\[
{\rm det}(v_1,\ldots,v_s)=\sum_{(j_1,\ldots,j_s)\in[q]^s}
\lambda_{1, j_1}\cdots\lambda_{s,j_s}\times
{\rm det}(w_{j_1},\ldots,w_{j_s})\ .
\] 
Assuming we have proved our inequality when the vectors are in $E_s$, we deduce
\begin{eqnarray*}
\left|{\rm det}(v_1,\ldots,v_s)\right| & \le & 
\sum_{(j_1,\ldots,j_s)\in[q]^s}
\lambda_{1, j_1}\cdots\lambda_{s,j_s}\times
\left|{\rm det}(w_{j_1},\ldots,w_{j_s})\right| \\
 & \le & \sum_{(j_1,\ldots,j_s)\in[q]^s}
\lambda_{1, j_1}\cdots\lambda_{s,j_s} \\
 & = & 1\ .
\end{eqnarray*}
We now prove the inequality when the vectors $v_1,\ldots,v_s$
belong to $E_s$, by induction on $s$. The case $s=1$ is trivial, so now assume $s\ge 2$. Let $R$ denote the matrix with rows given by the vectors $v_1,\ldots,v_s$. Permuting columns, or rows, or changing the sign of some rows does not affect the desired conclusion $|{\rm det}(R)|\le 1$, nor the hypotheses at our disposal.
If some rows are repeated (up to a sign) then ${\rm det}(R)=0$ and we are done, so we now assume there are no such repeats. If all the rows are of the form $e_i-e_j$, then $R(1,\ldots,1)^{\rm T}=0$, and we again have ${\rm det}(R)=0$. We are left with the case where there are $t\ge 1$ rows of the form $\pm e_i$. Since none of the corresponding indices $i$ are repeated, we can use the three types of operations mentioned above to bring $R$ into the block form
\[
R=\begin{pmatrix}
{\rm I}_t & 0 \\
S & T
\end{pmatrix}\ ,
\]
where ${\rm I}_t$ is the $t\times t$ identity matrix.
If $t=s$, then ${\rm det}(R)=1$, and we are done. Otherwise, we now have $1\le t\le s-1$, and we examine the $(s-t)\times(s-t)$ matrix $T$. If any of its rows is full of zeroes, then its determinant vanishes, and therefore ${\rm det}(R)=0$. If not, then its rows, which are of length $s-t$, are either of the form
\[
(0,\ldots,0,\pm 1,0,\ldots,0)\ ,
\]
or
\[
\pm(0,\ldots,0, -1, 0,\ldots,0,1,0,\ldots,0)\ .
\]
In other words, $T$ satisfies the same hypotheses as $R$ did, with respect to the size $s-t$ instead of $s$, and by induction $|{\rm det}(T)|\le 1$. Hence $|{\rm det}(R)|\le 1$.
\qed

Note that when proving the Perron-Frobenius Theorem, the use of the submultiplicative property of the matrix norm $\|P\|_{\infty,1}:=\max_{i}\sum_{j}|P_{ij}|$ immediately shows that the eigenvalues of a substochastic matrix are in the closed unit disk. Hence $|{\rm det}(P)|\le 1$. However, differences of substochastic matrices can have eigenvalues outside the unit disk, as shown by the simple example
\[
P=\begin{pmatrix}
1 & 0 \\
0 & 1
\end{pmatrix}
\ \ ,\ \ 
Q=\begin{pmatrix}
0 & 1 \\
1 & 0
\end{pmatrix}\ .
\]
Note that via the inverse of the linear isomorphism $(x_0,x_1,\ldots,x_s)\mapsto (x_1,\ldots,x_s)$ from $\{(x_0,x_1,\ldots,x_s)\in\mathbb{R}^{s+1}\ |\ x_0+x_1+\cdots+x_s=0\}$ to $\mathbb{R}^s$, the convex set $K_s$ can be seen as a full root polytope in type $A$ representation theory~\cite{Cho}. It is different from, yet related to the more common root polytope (for positive roots) introduced in~\cite{GelfandGP}.
Another proof of the determinantal inequality, in the case of vectors in $E_s$, is given in~\cite[Prop. 2.5]{Seashore}, using the notion of totally unimodular matrices (see, e.g.,~\cite[Ch. 19 to 21]{Schrijver}).
Finally, note that by multiplying the matrix $H$, on the left, by the inverse of the diagonal matrix in (\ref{diagonaleps}), before removing the row and column, and using Theorem \ref{substochthm} with $s=n-1$, we immediately obtain a proof of the inequality (\ref{minorineq}), which motivates the terminology of extremal pairs of transversals.

\section{The formula for the coefficients of resultants, by an induction on the dimension}
\label{formulasec}

We now tackle our main task of deriving a formula for the resultant. We assume $n\ge 2$, and the format defining the diagram $\mathbb{D}$ is $(n,n,(d_1,\ldots,d_n)^{\rm T})$ as in \S\ref{warmupsec}. The multidegree is the canonical one, namely, corresponding to $\delta=(\delta_1,\ldots,\delta_n)^{\rm T}$, where for all $i\in[n]$, $\delta_i:=\prod_{j\in[n]\backslash\{i\}}d_j$. The IFDZ $\mathscr{R}$, of format $(n,n,d)$ and multidegree $\delta^{\rm T}$, we will focus on is the resultant
\[
\mathscr{R}:={\rm Res}_{d_1,\ldots,d_n}\ .
\]
Taking Theorem \ref{balancedthm} into account, its expansion is given by
\[
\mathscr{R}=\sum_{A\in{\rm REF}} r_A\ \ua^{A}\ ,
\]
where
\[
{\rm REF}:=\{A\in{\rm EF}(\delta)\ |\ {\rm wt}(A)=d_1\cdots d_n(1,\ldots,1)\}
\]
refers to the set of resultant effective fillings, or effective fillings relevant to the computation of the resultant. Among those, we also define the ones which are {\it first-row-reduced} by
\[
{\rm FRR}:=\{A\in{\rm REF}\ |\ A_{1,d_1 e_1}=\delta_1\}\ .
\]
The following lemma is easy, yet important for the rest of the section.

\begin{lemma}\label{FRRproplem}
An effective first-row-reduced filling $A$ must satisfy $A_{1,\alpha}=0$, for all $\alpha\in\widetilde{\mathbb{D}}_1\backslash\{d_1 e_1\}$.
It must also satisfy $A_{i,\alpha}=0$ for all $i\in\{2,\ldots,n\}$, and all multiindex $\alpha=(\alpha_1,\ldots,\alpha_n)$ such that $\alpha_1>0$.
\end{lemma}

\noindent{\bf Proof:}
Indeed,
$\sum_{\alpha\in\widetilde{\mathbb{D}}_1}A_{1,\alpha}=\delta_1$, the nonnegativity of the entries of $A$, and the hypothesis $A_{1,d_1 e_1}=\delta_1$ immediately imply the vanishing of $A_{1,\alpha}$ for $\alpha\neq d_1 e_1$.
We also have
\[
{\rm wt}(A)=\delta_1 d_1 e_1+\sum_{i=2}^{n}\sum_{\alpha\in\widetilde{\mathbb{D}}_i} A_{i,\alpha}\ \alpha\ ,
\]
which, by taking first components, gives
\[
d_1\cdots d_n=d_1\cdots d_n+\sum_{i=2}^{n}\sum_{\alpha\in\widetilde{\mathbb{D}}_i} A_{i,\alpha}\ \alpha_1\ .
\]
As a result, $\forall i\in\{2,\ldots,n\}$, $\forall \alpha\in\widetilde{\mathbb{D}}_i$, $A_{i,\alpha}\ \alpha_1=0$, and the claim follows. 
\qed

We proceed as in the proof of Theorem \ref{balancedthm}, as far as setting up the resolvent equations. We take the same transversals $\mathbb{E}$ and $\mathbb{T}$ defined in (\ref{diagonaleps}) and (\ref{taupick}). Assume the filling $A$ is in ${\rm REF}\backslash{\rm FRR}$. We successively define $M,C,\beta$ using (\ref{resolvpick}). We also define the counter in the same way by $\mathcal{Z}(A):=\|\mathbf{1}_{\mathbb{D}\backslash(\mathbb{E}\cup\mathbb{T})}A\|$. The notion of $\wtA$ featuring in the resolvent equation for $A$ is defined in the same way as in \S\ref{warmupsec}. We will, however, slightly modify the definition of the strict partial order $\prec$, which will now be a relation on ${\rm REF}$. For $A,B\in {\rm REF}$, we will say $A\prec B$ iff
\[
(A\in{\rm FRR}\  {\rm and}\ B\notin{\rm FRR})\ {\rm or}\ (A,B\notin{\rm FRR}\ {\rm and}\ \mathcal{Z}(A)<\mathcal{Z}(B))\ .
\]
We then have the following analogue of Lemma \ref{keylem}.

\begin{lemma}\label{reskeylem}
Let $A\in{\rm REF}\backslash{\rm FRR}$ and $\wtA\in{\rm REF}$. If $\wtA\neq A$ and $\wtA$ features in the resolvent equation for $A$, then $\wtA\prec A$.
\end{lemma}

\noindent{\bf Proof:}
If the counter drops, i.e., $\mathcal{Z}(\wtA)<\mathcal{Z}(A)$, then, by our choice of definition, $\wtA\prec A$, regardless of whether $\wtA$ is first-row-reduced or not. So we now assume $\mathcal{Z}(\wtA)=\mathcal{Z}(A)$,
and repeat the same steps (using the same notation)  as in the proof of Theorem \ref{balancedthm}, therefore ariving at the same two cases to discuss. 
If $y_n=-1$, we again deduce $A_{1,d_1 e_1}=\delta_1$, i.e., $A\in{\rm FRR}$. However, this is excluded by hypothesis. So the second case is forced, namely, we must have $y_n=1$, but this again leads to $\wtA_{1,d_1 e_1}=\delta_1$. Hence, $\wtA\in{\rm FRR}$ and $\wtA\prec A$ by definition.
\qed

In the next proposition, we will rewrite the resolvent equation for a given $A\in{\rm REF}\backslash{\rm FRR}$ in a form suitable for iteration.
For that purpose, we define the transition matrix $\mathscr{T}$, as follows.
For any $A\in{\rm REF}\backslash{\rm FRR}$, and $\wtA\in{\rm REF}$, we let
\begin{eqnarray}
\mathscr{T}(A,\wtA) &:=&
\bbone\left\{
\mathbf{1}_{\mathbb{D}\backslash\mathbb{E}}\wtA
\le \mathbf{1}_{\mathbb{D}\backslash\mathbb{E}}A+P_{\mathbb{T}}(\mathbf{1}_{\mathbb{E}}A)
\right\} \nonumber \\
 & & \times (-1)^{\|\mathbf{1}_{\mathbb{E}}A\|+\|\mathbf{1}_{\mathbb{E}}\wtA\|+1}
\ 
\begin{bmatrix}
{\rm Row}(\mathbf{1}_{\mathbb{E}}\wtA) \\
\mathbf{1}_{\mathbb{D}\backslash\mathbb{E}}(A-\wtA)+P_{\mathbb{T}}(\mathbf{1}_{\mathbb{E}}A)
\end{bmatrix}\ .
\label{Tdefeq}
\end{eqnarray}

\begin{proposition}\label{Markovprop}
For all $A\in{\rm REF}\backslash{\rm FRR}$, we have
\[
r_A=\sum_{\substack{\wtA\in{\rm REF} \\ \wtA\prec A}}
\mathscr{T}(A,\wtA)\ r_{\wtA}\ .
\]
\end{proposition}

\noindent{\bf Proof:}
We consider the resolvent equation ${\rm Eq}(\mathbb{E},C,\beta)$ for $A$, as defined above, and write it in the $A$ version (\ref{longmaineqA}), using $\wtA$ as a summation index, instead of $A$ which is fixed. It is easy to see that the term where $\wtA=A$ contributes
$(-1)^{\|\mathbf{1}_{\mathbb{E}}A\|}\ r_A$,
because
\[
\begin{bmatrix}
{\rm Row}(\mathbf{1}_{\mathbb{E}}A) \\
C-\mathbf{1}_{\mathbb{D}\backslash\mathbb{E}}A
\end{bmatrix}
=
\begin{bmatrix}
{\rm Row}(\mathbf{1}_{\mathbb{E}}A) \\
M
\end{bmatrix}
=
\begin{bmatrix}
{\rm Row}(\mathbf{1}_{\mathbb{E}}A) \\
P_{\mathbb{T}}(\mathbf{1}_{\mathbb{E}}A)
\end{bmatrix}
=1\ ,
\]
by design.
For $\wtA\neq A$, the contribution is
\[
\bbone\left\{
\begin{array}{c}
\mathbf{1}_{\mathbb{D}\backslash\mathbb{E}}\wtA\le C \\
{\rm wt}(\wtA)={\rm wt}(C)-\beta+\delta^{\rm T}\varepsilon
\end{array}
\right\}
(-1)^{\|\mathbf{1}_{\mathbb{E}}\wtA\|}
\begin{bmatrix}
{\rm Row}(\mathbf{1}_{\mathbb{E}}\wtA) \\
C-\mathbf{1}_{\mathbb{D}\backslash\mathbb{E}}\wtA
\end{bmatrix}
r_{\wtA}\ .
\]
Since, by construction via (\ref{resolvpick}), ${\rm wt}(A)={\rm wt}(C)-\beta+\delta^{\rm T}\varepsilon$, the weight condition on $\wtA$ is ${\rm wt}(\wtA)={\rm wt}(A)=d_1\ldots d_n(1,\ldots,1)$, which is taken care of by summing $\wtA$ over ${\rm REF}$ instead of ${\rm EF}(\delta)$.
Since $C=\mathbf{1}_{\mathbb{D}\backslash\mathbb{E}}A+M$, and $M=P_{\mathbb{T}}(\mathbf{1}_{\mathbb{E}}A)$, we have
\[
C-\mathbf{1}_{\mathbb{D}\backslash\mathbb{E}}\wtA=
\mathbf{1}_{\mathbb{D}\backslash\mathbb{E}}(A-\wtA)+P_{\mathbb{T}}(\mathbf{1}_{\mathbb{E}}A)\ .
\]
Also note that $\wtA\neq A$ becomes the more precise statement $\wtA\prec A$, in light of Lemma \ref{reskeylem}. Easy algebra and clean up thus allows us to rewrite the resolvent equation in the form stated in the proposition.
\qed

By iterating Proposition \ref{Markovprop}, we immediately obtain the following corollary.
Before stating it, we define for any $A\in{\rm REF}$, and $B\in{\rm FRR}$, the new matrix element
\begin{equation}
\mathscr{C}(A,B)  := \bbone\{A\in{\rm FRR}, A=B\}\ + 
\bbone\{A\notin{\rm FRR} \}\times
\sum_{k\ge 0} \mathscr{T}^{k+1}(A,B)\ ,
\label{CABdefeq}
\end{equation}
where
\begin{eqnarray}
\mathscr{T}^{k+1}(A,B) & := &
\sum_{\wtA_1,\ldots,\wtA_k\in{\rm REF}\backslash{\rm FRR}}
\bbone\left\{
B\prec\wtA_k\prec\cdots\prec \wtA_1\prec A
\right\}\times \nonumber \\
 & & \qquad\qquad
\mathscr{T}(A,\wtA_1)
\mathscr{T}(\wtA_1,\wtA_2)\cdots
\mathscr{T}(\wtA_{k-1},\wtA_k)
\mathscr{T}(\wtA_k,B)\ .
\label{Tpowereq}
\end{eqnarray}
Note that for $k=0$, this just reduces to $\mathscr{T}^{1}(A,B)=\mathscr{T}(A,B) $.

\begin{corollary}
For any $A\in{\rm REF}$, we have
$r_A=\sum_{B\in{\rm FRR}}\mathscr{C}(A,B)\ r_B$.
\end{corollary}

Note that the sum over $k$ terminates, since one has the easy bound $k\le \mathcal{Z}(A)$, when $A\notin{\rm FRR}$, because of our definition of the strict partial order and the requirement that $\wtA_1,\ldots,\wtA_{k}$ are not in ${\rm FRR}$.
We reduced the computation of $r_A$'s from the general case $A\in{\rm REF}$ to that of $A$'s which are first-row-reduced. In order to proceed, we appeal to the identity (\ref{Laplaceeq}) from the introduction.
For all $i$, such that $2\le i\le n$, we define
\[
\widetilde{\mathbb{D}}_{i}^{(2)}:=
\{\alpha\in\mathbb{N}^n\ |\ \alpha_1=0\}\ ,
\]
where, now and later, we will use the notation $\alpha_j$ for the $j$-th component of a multiindex $\alpha\in\mathbb{N}^n$.
We then define the subdiagram
\[
\mathbb{D}^{(2)}:=\{(i,\alpha)\in\mathbb{D}\ |\ 2\le i\le n\ ,\ \alpha\in\widetilde{\mathbb{D}}_{i}^{(2)}\}\ .
\]
The precise definition of the forms $\overline{F}_{2}^{(2)},\ldots,\overline{F}_{n}^{(2)}$ appearing in (\ref{Laplaceeq}) is
\[
\overline{F}_{i}^{(2)}=\sum_{\alpha\in\widetilde{\mathbb{D}}_{i}^{(2)}}
a_{i,\alpha}\ \ux^{\alpha}\ ,
\]
for $2\le i\le n$.
These can be thought of as universal or generic homogeneous polynomials of degrees
$d_2,\ldots,d_n$, in the variables $x_2,\ldots,x_n$. Instead of renumbering them from $1$ to $n-1$, we prefer to keep their native numbering from $2$ to $n$ and to see their coefficients, which still are formal variables, collectively denoted $\ua_{\mathbb{D}^{(2)}}=(a_{i,\alpha})_{(i,\alpha)\in\mathbb{D}^{(2)}}$, as a subcollection
of the original variables $\ua_{\mathbb{D}}$.
Now, the expression
\[
{\rm Res}_{d_1,\ldots,d_n}\left(x_1^{d_1},\overline{F}_{2}^{(2)},\ldots,\overline{F}_{n}^{(2)}\right)
\]
is obtained from
\[
{\rm Res}_{d_1,\ldots,d_n}(F_1,\ldots,F_n)=\mathscr{R}=\sum_{A\in{\rm REF}}r_A
\ \ua^{A}\ ,
\]
by the following partial specialization.
We set
\[
a_{1,\alpha}:=\bbone\{\alpha=d_1 e_1\}\ ,
\]
for all $\alpha\in\widetilde{\mathbb{D}}_{1}$, and we also set
\[
a_{i,\alpha}:=0\ ,
\]
for all $i\in\{2,\ldots,n\}$, and multiindex $\alpha$ in $\widetilde{\mathbb{D}}_{i}\backslash\widetilde{\mathbb{D}}_{i}^{(2)}$, i.e., for which $\alpha_1>0$.
All the other $a_{i,\alpha}$ remain unspecialized formal variables.
From Lemma \ref{FRRproplem}, it is easy to see that
\[
{\rm Res}_{d_1,\ldots,d_n}\left(x_1^{d_1},\overline{F}_{2}^{(2)},\ldots,\overline{F}_{n}^{(2)}\right)=\sum_{A\in{\rm FRR}}r_A\ \ua^{\mathbf{1}_{\mathbb{D}^{(2)}}A}\ .
\]
On the other hand, we have
\begin{equation}
{\rm Res}_{d_2,\ldots,d_n}\left(\overline{F}_{2}^{(2)},\ldots,\overline{F}_{n}^{(2)}\right)
=:\sum_{B\in{\rm REF}^{(2)}}
r_{B}^{(2)}\ \ua^{B}\ ,
\label{smallreseq}
\end{equation}
where we wrote $=:$ in order to emphasize that this equation serves as definition for the coefficients $r_{B}^{(2)}$, rather than for the left-hand side. The latter is just the standard multidimensional resultant with generic or universal coefficients, except for numbering the forms and variables from $2$ to $n$, and seeing the variables $\ua_{\mathbb{D}^{(2)}}$ of the smaller resultant ${\rm Res}_{d_2,\ldots,d_n}$, as a subcollection of the original variables $\ua_{\mathbb{D}}$.
The set ${\rm REF}^{(2)}$ is the analogue, in this size $n-1$ relative setting, of the previous set ${\rm REF}$, in size $n$. In order to properly specify it, we first define
\[
\delta^{(2)}:=\begin{pmatrix}
\delta_{1}^{(2)} \\
\vdots \\
\delta_{n}^{(2)}
\end{pmatrix}\ ,
\]
where $\delta_{1}^{(2)}:=0$, and
\[
\delta_{i}^{(2)}:=\prod_{\substack{j=2 \\ j\neq i}}^{n}d_j \ ,
\]
for all $i$, such that $2\le i\le n$.
By definition, ${\rm REF}^{(2)}$ is the set of maps $B:\mathbb{D}\rightarrow\mathbb{B}$, such that ${\rm supp}(B)\subset\mathbb{D}^{(2)}$, ${\rm Row}(B)=\delta^{(2)}$, and
\[
{\rm wt}(B)=d_2\cdots d_n(0,1,\ldots,1)\ .
\]
As we will soon see in the iteration of the dimensional reduction embodied by Proposition \ref{dimredprop} below, it is convenient to avoid defining the smaller relative fillings $B$ as maps $\mathbb{D}^{(2)}\rightarrow\mathbb{N}$, and instead use the familiar trick of completing with zeros, so as to see $B$ as defined on the larger domain $\mathbb{D}$. 
Using (\ref{Laplaceeq}), (\ref{smallreseq}), and expanding the $d_1$-th power, we obtain
\[
{\rm Res}_{d_1,\ldots,d_n}\left(x_1^{d_1},\overline{F}_{2}^{(2)},\ldots,\overline{F}_{n}^{(2)}\right)
=\sum_{B_1,\ldots,B_{d_1}\in {\rm REF}^{(2)}}
r_{B_1}^{(2)}\cdots r_{B_{d_1}}^{(2)}\ \ua^{B_1+\cdots+B_{d_1}}\ .
\]
By extracting the coefficient of $\ua^{\mathbf{1}_{\mathbb{D}^{(2)}}A}$, we see that, for all
$A\in{\rm FRR}$,
\begin{equation}
r_A=\sum_{B_1,\ldots,B_{d_1}\in {\rm REF}^{(2)}}
\bbone\left\{
B_1+\cdots+B_{d_1}=\mathbf{1}_{\mathbb{D}^{(2)}}A
\right\}\ 
r_{B_1}^{(2)}\cdots r_{B_{d_1}}^{(2)}\ .
\label{FRRformulaeq}
\end{equation}
In order to consolidate our notation, and write more compact formulas, we define
\[
\mathscr{W}^{(1)}(A;B_1,\ldots,B_{d_1}):=
\sum_{B\in{\rm FRR}}
\bbone\left\{
B_1+\cdots+B_{d_1}=\mathbf{1}_{\mathbb{D}^{(2)}}B
\right\}\ 
\mathscr{C}(A,B)\ ,
\]
for all $A\in{\rm REF}$ and all $B_1,\ldots,B_{d_1}\in{\rm REF}^{(2)}$.
We have established the following proposition which gives a dimensional reduction, from $n$ to $n-1$, for the computation of the coefficients of the multidimensional resultant.

\begin{proposition}\label{dimredprop}
For all $A\in{\rm REF}$, we have
\[
r_A=\sum_{B_1,\ldots,B_{d_1}\in {\rm REF}^{(2)}}
\mathscr{W}^{(1)}(A;B_1,\ldots,B_{d_1})\ 
r_{B_1}^{(2)}\cdots r_{B_{d_1}}^{(2)}\ .
\] 
\end{proposition}

Iterating this identity, from dimension $n$ all the way to dimension $1$, presents no difficulty, except notation management. We hope that the following choices of definitions and notations will appear natural to the reader, as well as consistent with the choices made so far.
Let $\ell\in[n]$, and define
\[
\widetilde{\mathbb{D}}_{i}^{(\ell)}:=
\{\alpha\in\mathbb{N}^n\ |\ \alpha_1=0,\ldots,\ \alpha_{\ell-1}=0,\ |\alpha|=d_i\}\ ,
\]
for all $i$ such that $\ell\le i\le n$.
We also define the subdiagram
\[
\mathbb{D}^{(\ell)}:=\{(i,\alpha)\in\mathbb{D}\ |\ \ell\le i\le n,\ \alpha\in\widetilde{\mathbb{D}}_{i}^{(\ell)}\}\ ,
\]
as well as the column vector
\[
\delta^{(\ell)}:=\begin{pmatrix}
\delta_{1}^{(\ell)} \\
\vdots \\
\delta_{n}^{(\ell)}
\end{pmatrix}\ ,
\]
where
\[
\delta_{i}^{(\ell)}:=
\left\{
\begin{array}{cl}
0 & {\rm if}\ 1\le i<\ell\ , \\
\prod\limits_{\substack{j=\ell \\ j\neq i}}^{n}d_j & {\rm if}\ \ell\le i\le n \ .
\end{array}\right.
\]
We define the relative notion of resultant effective fillings for the subcollection of $x$ variables and homogeneous forms, numbered from $\ell$ to $n$, as follows.
We let ${\rm REF}^{(\ell)}$ denote the set of maps $B:\mathbb{D}\rightarrow\mathbb{N}$, such that ${\rm supp}(B)\subset\mathbb{D}^{(\ell)}$, ${\rm Row}(B)=\delta^{(\ell)}$, and
\[
{\rm wt}(B)=d_{\ell}\cdots d_n(0,\ldots,0,1,\ldots,1)\ ,
\]
where the first $1$ appears in the $\ell$-th position.
Next, we define the relative notion of first-row-reduced fillings, namely, we let
\[
{\rm FRR}^{(\ell)}:=\{B\in{\rm REF}^{(\ell)}\ |\ B_{\ell,d_{\ell} e_{\ell}}=\delta_{\ell}^{(\ell)}\}\ .
\]
We also define the coefficients $r_{B}^{(\ell)}$ of the resultant of size $n-\ell+1$, by
\[
{\rm Res}_{d_{\ell},\ldots,d_n}\left(\overline{F}_{\ell}^{(\ell)},\ldots,\overline{F}_{n}^{(\ell)}\right)
=:\sum_{B\in{\rm REF}^{(\ell)}}
r_{B}^{(\ell)}\ \ua^{B}
\]
where the
\[
\overline{F}_{i}^{(\ell)}:=\sum_{\alpha\in\widetilde{\mathbb{D}}_{i}^{(\ell)}} a_{i,\alpha}\ \ux^{\alpha}\ ,
\]
for $\ell\le i\le n$,
are the generic or universal homogeneous polynomials of degrees $d_{\ell},\ldots,d_n$ in the variables $x_{\ell},\ldots,x_n$.

In the particular case $\ell=1$, the objects $\widetilde{\mathbb{D}}_{i}^{(\ell)}$,
$\mathbb{D}^{(\ell)}$, $\delta^{(\ell)}$, ${\rm REF}^{(\ell)}$, ${\rm FRR}^{(\ell)}$, and $\overline{F}_{i}^{(\ell)}$
respectively coincide with what we previously simply denoted by
$\widetilde{\mathbb{D}}_{i}$,
$\mathbb{D}$, $\delta$, ${\rm REF}$, ${\rm FRR}$, and $F_i$.
The previously introduced $\widetilde{\mathbb{D}}_{i}^{(2)}$,
$\mathbb{D}^{(2)}$, $\delta^{(2)}$, ${\rm REF}^{(2)}$, and $\overline{F}_{i}^{(2)}$ are subsumed under our new definitions, as the $\ell=2$ particular case.

It is important to consider the rather degenerate yet meaningful case $\ell=n$, since it is the initial condition for our iterative procedure. The resultant of a single homogeneous polynomial in one variable still makes sense (like the determinant of a $1\times 1$ matrix) and, with self-explanatory notations, is given, in our context by
\begin{equation}
{\rm Res}_{d_n}\left(\overline{F}_{n}^{(n)}\right)=
{\rm Res}_{d_n}\left(a_{n,d_n e_n}x_n^{d_n}\right)=a_{n,d_n e_n}\ .
\label{trivresn}
\end{equation}
Note that
\[
\delta^{(n)}:=\begin{pmatrix}
0 \\
\vdots \\
0 \\
1
\end{pmatrix}\ ,
\]
and all the previous sets become singletons:
\begin{eqnarray*}
\widetilde{\mathbb{D}}_{n}^{(n)} & = & \{d_n e_n\} \ , \\
\mathbb{D}^{(n)} & = & \{(n,d_n e_n)\} \ , \\
{\rm REF}^{(n)} & = & \{\mathbf{1}_{\mathbb{D}^{(n)}}\} \ .
\end{eqnarray*}
There is only one coefficient to consider for the $\ell=n$ case of our relative resultants, and (\ref{trivresn}) provides us with our {\it initial condition}
\begin{equation}
r_{\mathbf{1}_{\mathbb{D}^{(n)}}}^{(n)}=1\ .
\label{initcondeq}
\end{equation}
When $1\le \ell\le n-1$, we define relative versions of the transversals $\mathbb{E}$, $\mathbb{T}$ from (\ref{diagonaleps}) and (\ref{taupick}).
Namely, we let
\[
\mathbb{E}^{(\ell)}:=\{(i,d_i e_i)\ |\ \ell\le i\le n\}\ ,
\]
and
\[
\mathbb{T}^{(\ell)}:=\{(i,(d_i-1)e_{\ell}+e_{i+1})\ |\ \ell\le i<n\}\ \bigcup\ 
\{(n,d_n e_{\ell})\}\ .
\]
We also define a relative push operator to $\mathbb{T}^{(\ell)}$ by letting, for effective fillings $A$ such that ${\rm supp}(A)\subset\mathbb{D}^{(\ell)}$,
\[
\left[P_{\mathbb{T}^{(\ell)}}^{(\ell)}(A)\right]_{i,\alpha}:=
{\rm Row}_i(A)\times \bbone\{(i,\alpha)\in\mathbb{T}^{(\ell)}\}\ ,
\]
for all $(i,\alpha)\in\mathbb{D}$.
We define a relative counter
\[
\mathcal{Z}^{(\ell)}(A):=\|\mathbf{1}_{\mathbb{D}^{(\ell)}\backslash(\mathbb{E}^{(\ell)}\cup \mathbb{T}^{(\ell)})}A\|\ ,
\]
for $A$'s in ${\rm REF}^{(\ell)}$.
Next, we define the relative strict partial order $A\prec^{(\ell)} B$, for $A,B\in{\rm REF}^{(\ell)}$, by the statement
\[
(A\in{\rm FRR}^{(\ell)}\  {\rm and}\ B\notin{\rm FRR}^{(\ell)})\ {\rm or}\ (A,B\notin{\rm FRR}^{(\ell)}\ {\rm and}\ \mathcal{Z}^{(\ell)}(A)<\mathcal{Z}^{(\ell)}(B))\ .
\]
For any $A\in{\rm REF}^{(\ell)}\backslash{\rm FRR}^{(\ell)}$, and $\wtA\in{\rm REF}^{(\ell)}$, we let
\begin{eqnarray*}
\mathscr{T}^{(\ell)}(A,\wtA) & := &
\bbone\left\{
\mathbf{1}_{\mathbb{D}^{(\ell)}\backslash\mathbb{E}^{(\ell)}}\wtA
\le \mathbf{1}_{\mathbb{D}^{(\ell)}\backslash\mathbb{E}^{(\ell)}}A+
P_{\mathbb{T}^{(\ell)}}^{(\ell)}(\mathbf{1}_{\mathbb{E}^{(\ell)}}A)
\right\} \times \\
 & &
\ (-1)^{\|\mathbf{1}_{\mathbb{E}^{(\ell)}}A\|+\|\mathbf{1}_{\mathbb{E}^{(\ell)}}\wtA\|+1}
\ 
\begin{bmatrix}
{\rm Row}(\mathbf{1}_{\mathbb{E}^{(\ell)}}\wtA) \\
\mathbf{1}_{\mathbb{D}^{(\ell)}\backslash\mathbb{E}^{(\ell)}}(A-\wtA)+P_{\mathbb{T}^{(\ell)}}^{(\ell)}(\mathbf{1}_{\mathbb{E}^{(\ell)}}A)
\end{bmatrix}\ .
\end{eqnarray*}
For any $A\in{\rm REF}^{(\ell)}$, and $B\in{\rm FRR}^{(\ell)}$, we define
\begin{eqnarray*}
\mathscr{C}^{(\ell)}(A,B) & := &\bbone\{A\in{\rm FRR}^{(\ell)}, A=B\}\ + \ 
\bbone\{A\notin{\rm FRR}^{(\ell)} \}\times \\
 & & \qquad
\sum_{k\ge 0}\sum_{\wtA_1,\ldots,\wtA_k\in{\rm REF}^{(\ell)}\backslash{\rm FRR}^{(\ell)}}
\bbone\left\{
B\prec^{(\ell)}\wtA_k\prec^{(\ell)}\cdots\prec^{(\ell)} \wtA_1\prec^{(\ell)} A
\right\}\times \\
 & & \qquad\qquad
\mathscr{T}^{(\ell)}(A,\wtA_1)
\mathscr{T}^{(\ell)}(\wtA_1,\wtA_2)\cdots
\mathscr{T}^{(\ell)}(\wtA_{k-1},\wtA_k)
\mathscr{T}^{(\ell)}(\wtA_k,B)\ .
\end{eqnarray*}

\smallskip
Finally, for all $A\in{\rm REF}^{(\ell)}$ and all $B_1,\ldots,B_{d_l}\in{\rm REF}^{(\ell+1)}$, we define
\[
\mathscr{W}^{(\ell)}(A;B_1,\ldots,B_{d_{\ell}}):=
\sum_{B\in{\rm FRR}^{(\ell)}}
\bbone\left\{
B_1+\cdots+B_{d_{\ell}}=\mathbf{1}_{\mathbb{D}^{(\ell+1)}}B
\right\}\ 
\mathscr{C}^{(\ell)}(A,B)\ .
\]
From the same arguments which lead to Proposition \ref{dimredprop} in the context of the full resultant ${\rm Res}_{d_1,\ldots,d_n}(F_1,\ldots,F_n)$, one can deduce the following result about the coefficients of smaller intermediate resultants ${\rm Res}_{d_{\ell},\ldots,d_n}\left(\overline{F}_{\ell}^{(\ell)},\ldots,\overline{F}_{n}^{(\ell)}\right)$.

\begin{proposition}\label{reldimredprop}
For all $\ell$, such that $1\le \ell\le n-1$, and
for all $A\in{\rm REF}^{(\ell)}$, we have
\[
r_A^{(\ell)}=\sum_{B_1,\ldots,B_{d_{\ell}}\in {\rm REF}^{(\ell+1)}}
\mathscr{W}^{(\ell)}(A;B_1,\ldots,B_{d_{\ell}})\ 
r_{B_1}^{(\ell+1)}\cdots r_{B_{d_{\ell}}}^{(\ell+1)}\ .
\] 
\end{proposition}

In order to write an induction-free and fully explicit formula for the coefficients of the resultant ${\rm Res}_{d_1,\ldots,d_n}(F_1,\ldots,F_n)$, the last ingredient we need is a rooted tree graph. The latter has a set of vertices decomposed into $n$ generations ${\rm Word}^{(\ell)}$, $1\le \ell\le n$. The first generation has only one vertex, the root of the tree. It has $d_1$ vertices attached to it belonging to the second generation. Each vertex of the latter is attached to $d_2$ vertices in the third generation, etc.
For the sake of precision, we will write ${\rm Word}^{(\ell)}:=[d_1]\times\cdots\times[d_{\ell-1}]$, and see an element or vertex $\mathsf{w}\in{\rm Word}^{(\ell)}$ as a `word' $\mathsf{w}=L_1L_2\cdots L_{\ell-1}$ made of a succession of $\ell-1$ `letters' $L_j$, $1\le j\le \ell-1$, such that $L_j\in[d_j]$.
In particular, ${\rm Word}^{(1)}=\{\varnothing\}$, where $\varnothing$ is the empty word and encodes the root vertex.

With all these preparatory steps done, we can state our main theorem, which immediately results from iterating Proposition \ref{reldimredprop}, until one hits the initial condition (\ref{initcondeq}).

\begin{theorem}\label{maintheorem}
For any $n\ge 1$, and $d_1,\ldots,d_n\ge 1$, the multidimentional resultant of the corresponding format is given by
\[
{\rm Res}_{d_1,\ldots,d_n}=\sum_{A\in{\rm REF}}r_A\ \ua^A\ ,
\] 
where, for all $A\in{\rm REF}$,
\[
r_A=\sum_{\mathscr{B}}\prod_{l=1}^{n-1}\left(
\prod_{\mathsf{w}\in{\rm Word}^{(\ell)}}
\mathscr{W}^{(\ell)}\left(
B_{\mathsf{w}}^{(\ell)}; 
B_{\mathsf{w}1}^{(\ell+1)},
B_{\mathsf{w}2}^{(\ell+1)},\ldots,
B_{\mathsf{w}d_{\ell}}^{(\ell+1)}
\right)
\right)\ .
\]
In the above, we used the concatenation notation where, for $\mathsf{w}\in{\rm Word}^{(\ell)}$, and $j\in[d_{\ell}]$, $\mathsf{w}j$ denotes the word in ${\rm Word}^{(\ell+1)}$ obtained from $\mathsf{w}$ by adding the letter $j$ at the end.
The sum is over all collections
\[
\mathscr{B}=(B_{\mathsf{w}}^{(\ell)})_{\substack{1\le \ell\le n\ \  \\
\mathsf{w}\in{\rm Word}^{(\ell)} }} \ \ \ \ ,
\]
such that, $\forall \ell\in[n]$, $\forall\mathsf{w}\in{\rm Word}^{(\ell)}$, $B_{\mathsf{w}}^{(\ell)}\in{\rm REF}^{(\ell)}$,
and which satisfy the constraint $B_{\varnothing}^{(1)}=A$.
\end{theorem}

Note that by unpacking the definitions of $\mathscr{W}^{(\ell)}$, etc. all the way to (\ref{supernomial}), it is easy to verify the claim made in the introduction about realizing the resultant coefficients as an alternating sum of products of multinomial coefficients. We thus obtained a new proof of the known fact that $r_A\in\mathbb{Z}$, for all $A\in{\rm REF}$.

\section{The linear case, or rediscovering determinants}
\label{linearsec}

In this section, we look at how the discussion from \S\ref{formulasec} specializes to the linear case, where $d_1=\cdots=d_n=1$. We now have $\mathbb{D}=\{(i,e_j)\ |\ i,j\in[n]\}$, and we therefore identify it with the square grid $[n]\times[n]$. We will also simplify our notation by writing $A_{i,e_j}=:A_{ij}$, thus seeing a filling $A$ as a familiar $n\times n$ matrix with nonnegative integer entries.
From ${\rm Row}(A)=(1,\ldots,1)^{\rm T}$, and ${\rm wt}(A)=(1,\ldots,1)$, we see that the entries of $A$ must be zeros or ones, and there is exactly a single one, in each row, and in each column. In other words, $A$ is a permutation matrix. We define the correspondence between such a matrix $A\in{\rm REF}$, and a permutation $\sigma$ in the symmetric group $\mathfrak{S}_n$, by letting
\[
A_{ij}=\bbone\{j=\sigma(i)\}\ ,
\]
for all $i,j\in[n]$.
We will take advantage of this one-to-one correspondence by rewriting our coefficients $r_{\sigma}$, instead of $r_A$.
In this section, we will see how to recover, from the discussion in \S\ref{formulasec}, the familiar fact that $r_{\sigma}$ is equal to the sign of the permutation $\sigma$.
The transversal $\mathbb{E}$ corresponds to the matrix $\varepsilon={\rm I}_n$, the identity matrix, whereas $\mathbb{T}$ corresponds to the matrix
\[
\tau=\begin{pmatrix}
0 & 1 & 0 & \cdots & 0 \\
0 & 0 & 1 & \ddots & \vdots \\
\vdots & \vdots & \ddots & \ddots & 0 \\
0 & 0 & \cdots & 0 & 1 \\
1 & 0 & \cdots & 0 & 0
\end{pmatrix}\ .
\]
Note that $A\in{\rm FRR}$ is equivalent to $\sigma(1)=1$.
We will denote by $f(\sigma)$ the number of fixed points of the permutation $\sigma$. It is easy to see that $f(\sigma)=\|\mathbf{1}_{\mathbb{E}}A\|$.
Now consider a matrix $A$ which is not first-row-reduced, i.e., such that $1$ is not fixed by the permutation $\sigma$.
The condition
\[
\mathbf{1}_{\mathbb{D}\backslash\mathbb{E}}\wtA
\le \mathbf{1}_{\mathbb{D}\backslash\mathbb{E}}A+P_{\mathbb{T}}(\mathbf{1}_{\mathbb{E}}A)\ ,
\]
in the definition of $\mathscr{T}(A,\wtA)$ from Proposition \ref{Markovprop}, now becomes
\begin{equation}
\wtA_{ij}\le A_{ij}+A_{ii}\ \bbone\{j=\rho(i)\}\ ,
\label{matrixpermcond}
\end{equation}
for all $i,j\in[n]$, such that $i\neq j$.
Here, $\rho$ is the maximal cycle $(12\cdots n)$ which sends $a\in[n]$ to $a+1$ modulo $n$, as seen from the matrix $\tau$.
Throughout the next few lemmas, $\wtA$ or rather the associated permutation $\wts$ is fixed and assumed to satisfy the conditions (\ref{matrixpermcond}).
We will also use $a$ to denote an element of $[n]$, rather than a formal variable $a_{i,\alpha}$, as in the previous sections.

\begin{lemma}\label{linlem1}
If $a\in[n]$ is such that $\sigma(a)\neq a$, and $\wts(a)\neq a$, then $\wts(a)=\sigma(a)$.
\end{lemma}

\noindent{\bf Proof:}
We apply (\ref{matrixpermcond}) with $i=a$ and $j=\wts(a)$, since the hypothesis provides for the requirement $i\neq j$.
We note that
\[
\wtA_{ij}=\bbone\{j=\wts(i)\}=\bbone\{\wts(a)=\wts(a)\}=1\ ,
\]
\[
A_{ij}=\bbone\{j=\sigma(i)\}=\bbone\{\wts(a)=\sigma(a)\}\ ,
\]
\[
A_{ii}=\bbone\{i=\sigma(i)\}=\bbone\{a=\sigma(a)\}=0\ .
\]
After simplification, (\ref{matrixpermcond}) reduces to
\[
1\le \bbone\{\wts(a)=\sigma(a)\}\ ,
\]
and therefore $\wts(a)=\sigma(a)$ must hold.
\qed

\begin{lemma}\label{linlem2}
Let $a\in[n]$, and suppose $\sigma(a)\neq a$, as well as $\wts(a)\neq a$, then
$\forall k\ge 0$, $\wts\left(\sigma^k(a)\right)=\sigma^{k+1}(a)$.
\end{lemma}

\noindent{\bf Proof:}
We proceed by induction on $k$. Lemma \ref{linlem1} establishes the $k=0$ base case.
Now suppose $k\ge 0$, and the induction hypothesis $\wts\left(\sigma^k(a)\right)=\sigma^{k+1}(a)$ holds.
Since $\sigma(a)\neq a$, and $\sigma^{k+1}$ is a bijection, we have $\sigma^{k+2}(a)\neq\sigma^{k+1}(a)$. Since $\wts\circ\sigma^k$ is bijective, we must also have
$\wts\left(\sigma^{k+1}(a)\right)\neq \wts\left(\sigma^{k}(a)\right)$, which can be rephrased as $\wts\left(\sigma^{k+1}(a)\right)\neq \sigma^{k+1}(a)$, via the induction hypothesis. Thus $\sigma^{k+1}(a)$, instead of $a$, satisfies the hypotheses of Lemma \ref{linlem1}, and we can conclude that $\wts\left(\sigma^{k+1}(a)\right)=\sigma^{k+2}(a)$, as desired.
\qed

\begin{lemma}\label{linlem3}
Let $a\in[n]$, and suppose $\sigma(a)\neq a$, as well as $\wts(a)= a$, then
$\forall k\ge 0$, $\wts\left(\sigma^k(a)\right)=\sigma^{k}(a)$.
\end{lemma}

\noindent{\bf Proof:}
Assume, for the sake of contradiction, that $\wts\left(\sigma^k(a)\right)\neq\sigma^{k}(a)$, for some $k\ge 0$.
Since $\sigma(a)\neq a$, and $\sigma^{k}$ is a bijection, we have $\sigma^{k+1}(a)\neq\sigma^{k}(a)$. 
Hence $\sigma^{k}(a)$, instead of $a$, satisfies the hypotheses of Lemma \ref{linlem2}. We therefore obtain $\forall\ell\ge 0$,
\[
\wts\left(\sigma^{\ell}\left(\sigma^k(a)\right)\right)=
\sigma^{\ell+1}\left(\sigma^k(a)\right)\ ,
\]
i.e.,
\[
\wts\left(\sigma^{k+\ell}(a)\right)=\sigma\left(\sigma^{k+\ell}(a)\right)\ .
\]
In other words, the permutations $\wts$, $\sigma$ coincide on the set $\{\sigma^{k+\ell}(a)\ |\ \ell\ge 0\}$ which is the support of the cycle of $\sigma$ which contains $a$. In particular, we must have $\wts(a)=\sigma(a)$ which contradicts the two hypotheses of the current lemma. As a result, the assumption made at the beginning must be rejected, and we must have $\wts\left(\sigma^k(a)\right)=\sigma^{k}(a)$, for all $k\ge 0$.
\qed

\begin{lemma}\label{linlem4}
Let $a\in[n]$, and suppose $\sigma(a)\neq a$.
Let $S=\{\sigma^{k}(a)\ |\ k\ge 0\}$ denote the support of the $\sigma$ cycle containing $a$. Then, we either have
\[
\forall b\in S,\ \ \wts(b)=\sigma(b)\ ,
\]
or
\[
\forall b\in S, \ \ \wts(b)=b\ .
\]
\end{lemma}

\noindent{\bf Proof:}
If $\wts(a)\neq a$, then Lemma \ref{linlem2} implies the first clause of the alternative.
Otherwise, $\wts(a)=a$, and Lemma \ref{linlem3} implies the second clause of the alternative.
\qed

\smallskip
We now note that $[n]$ is the disjoint union of sets $S_1,\ldots,S_q,F$, where the $S_i$, with $|S_i|\ge 2$, are the supports of the nontrivial cycles of $\sigma$, and $F$ is the set of fixed points of $\sigma$.
A consequence of Lemma \ref{linlem4} is that each subset $S_i$ is stable by $\wts$, i.e., satisfies $\wts(S_i)\subset S_i$. Since $\wts$ is a bijection, we get the more precise statement $\wts(S_i)=S_i$, for all $i$, such that $1\le i\le q$. Again because of the bijective property, we also deduce that $F$ is stable and, in fact satisfies $\wts(F)=F$.

\begin{lemma}\label{linlem5}
Let $a\in F$, and suppose that $\wts(a)\neq a$. Then, we have $\wts(a)=\rho(a)$. 
\end{lemma}

\noindent{\bf Proof:}
We apply (\ref{matrixpermcond}) with $i=a$, and $j=\wts(a)$.
We note that
\[
\wtA_{ij}=\bbone\{\wts(a)=\wts(a)\}=1\ ,
\]
trivially, while
\[
A_{ij}=\bbone\{\wts(a)=\sigma(a)\}=0\ ,
\]
because $a\in F$, i.e., $\sigma(a)=a$, and $\wts(a)\neq a$, by hypothesis.
Finally, we also have
\[
A_{ii}=\bbone\{a=\sigma(a)\}=1\ .
\]
After simplification, (\ref{matrixpermcond}) becomes
\[
1\le\bbone\{j=\rho(i)\}\ ,
\]
which implies the claim $\wts(a)=\rho(a)$.
\qed

\begin{lemma}\label{linlem6}
For all $a\in F$, we have $\wts(a)=a$.
\end{lemma}

\noindent{\bf Proof:}
We argue by contradiction, and assume that there exists some $a\in F$, such that
$\wts(a)\neq a$. Our first task is to show, by induction, that for all $k\ge 0$, the induction hypothesis
\[
\wts\left(\rho^k(a)\right)\neq \rho^k(a)\ \ ,\ {\rm and}\ \ \rho^k(a)\in F
\]
is satisfied.
The base case $k=0$ follows from the assumption of the reasoning by contradiction, and the hypothesis of the lemma.
Suppose the statement is true for $k\ge 0$. By Lemma \ref{linlem5}, applied to $\rho^k(a)$ instead of $a$, we deduce 
\begin{equation}
\wts\left(\rho^k(a)\right)=\rho^{k+1}(a)\ .
\label{linlem6eq}
\end{equation}
This implies $\rho^{k+1}(a)\in F$, because $F$ is stable by $\wts$, and $\rho^{k}(a)\in F$, by the induction hypothesis. Since the latter also includes $\wts\left(\rho^k(a)\right)\neq \rho^k(a)$, and $\wts$ is bijective, we obtain
\[
\wts^{2}\left(\rho^k(a)\right)\neq\wts\left(\rho^k(a)\right)\ ,
\]
which, using (\ref{linlem6eq}) on both sides, can be rewritten as
\[
\wts\left(\rho^{k+1}(a)\right)\neq\rho^{k+1}(a).
\]
This establishes the induction hypothesis for $k+1$, and completes the proof by induction.

For $k\ge 0$, we showed that $\rho^k(a)$, instead of $a$, satisfies the hypotheses of Lemma \ref{linlem5}, and therefore $\wts\left(\rho^k(a)\right)=\rho^{k+1}(a)$. This means that $\wts$ and $\rho$ coincide on the support of the cycle of $\rho$ containing $a$, i.e., all of $[n]$, because $\rho$ acts transitively. Hence $\wts=\rho$.
However, $a\in F$ and $F$ is stable by $\wts$ and therefore by $\rho$. This forces $F=[n]$, and $q=0$, which contradicts the standing hypothesis that $A\in{\rm REF}\backslash{\rm FRR}$, i.e., $\sigma(1)\neq 1$.
Our initial assumption $\wts(a)\neq a$ must be rejected. We showed that, for all $a\in F$, $\wts(a)=a$.
\qed

The previous series of lemmas can be summarized, using the matrix or filling identification with permutations, and the slight abuse of notation which goes with it, by the following proposition.

\begin{proposition}
If $\sigma\in{\rm REF}\backslash{\rm FRR}$, and the permutation $\wts\neq\sigma$ satisfies the conditions (\ref{matrixpermcond}) for all $i\neq j$, then one can obtain $\wts$ by deleting at least one nontrivial cycle from the cycle factorization of $\sigma$.
\end{proposition}

We can make the observation that
if $\wts$ is obtained from $\sigma\in{\rm REF}\backslash{\rm FRR}$ by cycle deletion as in the statement of the last proposition, then
\begin{equation}
\forall i\in[n],\ \ \left(\wts(i)=i\ \ {\rm or}\ \ \wts(i)=\sigma(i)\right)\ .
\label{observationeq}
\end{equation}

It is easy to see that the converse of the previous proposition is also true.
If $j\neq\wts(i)$, then $\wtA_{ij}=0$, and the inequality (\ref{matrixpermcond}) holds, regardless of the value of the right-hand side.
Now suppose $j=\wts(i)$. If $j=\sigma(i)$, then $A_{ij}=1$, and the inequality (\ref{matrixpermcond}) holds, regardless of the value of the left-hand side.
Finally, if $j\neq\sigma(i)$, then $\wts(i)\neq\sigma(i)$ which implies $\wts(i)=i$, by the observation (\ref{observationeq}).
This entails $j=i$, which is excluded from consideration, as far the conditions (\ref{matrixpermcond}) are concerned.
We have thus proved that, for all $\sigma\in{\rm REF}\backslash{\rm FRR}$,
\begin{equation}
r_{\sigma}=\sum_{\wts\in\mathfrak{S}_n}\mathscr{T}(\sigma,\wts) \ r_{\wts}\ ,
\label{perminducteq}
\end{equation}
where $\mathscr{T}(\sigma,\wts)$ is equal to the sign factor $(-1)^{f(\sigma)+f(\wts)+1}$ times the indicator of the condition that $\wts$ is obtained by deleting at least one nontrivial cycle from the cycle factorization of $\sigma$.

Using (\ref{perminducteq}), it is easy to prove that $r_{\sigma}$ is the sign $(-1)^{\sum_{i=1}^{q}(|S_i|-1)}$ of the permutation $\sigma$, e.g., by induction on the number of nontrivial cycles. Indeed, one can rewrite (\ref{perminducteq}) as
\begin{equation}
r[\sigma]=\sum_{\substack{\omega\in\{0,1\}^{[q]} \\ \omega\neq(0,\ldots,0)}}
(-1)^{1+\sum_{i=1}^{q}\omega_i|S_i|}
\ \ r\left[\prod_{i\in[q] , \omega_i=0}\rho_i\right]\ .
\label{perminduct2}
\end{equation}
In the above equation, we wrote the factorization of $\sigma$ into nontrivial cycles as $\sigma=\rho_1\circ\cdots\circ\rho_q$, where for all $i\in[q]$, $\rho_i$ has support $S_i$. For better readability, we changed the notation $r_{\sigma}$ to $r[\sigma]$.
The binary sequence $\omega=(\omega_1,\ldots,\omega_q)$ specifies the cycle deletion, namely, the cycle $\rho_i$ is deleted if $\omega_i=1$, and it is kept if $\omega_i=0$. Note the change in the number of fixed points is
\[
f(\wts)-f(\sigma)=\sum_{i=1}^{q}\omega_i|S_i|\ ,
\]
which explains the sign in (\ref{perminduct2}).
From the induction hypothesis, one can write the right-hand side of (\ref{perminduct2}) as
\begin{eqnarray*}
{\rm RHS} & = & \sum_{\substack{\omega\in\{0,1\}^{[q]} \\ \omega\neq(0,\ldots,0)}}
(-1)^{1+\sum_{i=1}^{q}\omega_i|S_i|}
\ \times\ (-1)^{\sum_{i=1}^{q}(1-\omega_i)(|S_i|-1)} \\
 & = & (-1)^{1+\sum_{i=1}^{q}(|S_i|-1)}
\times
 \sum_{\substack{\omega\in\{0,1\}^{[q]} \\ \omega\neq(0,\ldots,0)}}
\prod_{i=1}^{q}(-1)^{\omega_i} \\
 & = & (-1)^{1+\sum_{i=1}^{q}(|S_i|-1)}
\times\ 
\left[(1-1)^q-1\right] \\
 & = & (-1)^{\sum_{i=1}^{q}(|S_i|-1)}\ ,
\end{eqnarray*}
which propagates the induction hypothesis.

\begin{remark}
In the previous Remark \ref{singlexErem}, we mentioned that it is possible to solve for the coefficients $r_{\sigma}$, in the linear case, using a transversal $\mathbb{E}$ such as $\{(1,x_1^{d_1}),\ldots,(n,x_1^{d_n})\}$. We leave it to the diligent reader to check that the corresponding equations ${\rm Eq}(\mathbb{E},C,\beta)$ imply that, as one crosses an edge of the so-called star graph, the coefficient $r_{\sigma}$ picks up a minus sign. The star graph~\cite{AkersHK} is the Cayley graph of the symmetric group $\mathfrak{S}_n$ for the set of generators given by the transpositions $(1,i)$, with $2\le i\le n$. The key fact needed for this exercise, is that the prefix exchange distance (the graph distance on the star graph) to the identity permutation changes by $\pm 1$ as one moves along an edge (see, e.g.,~\cite[Lem. 1]{DayT}). This shows the graph is bipartite, and can be used to provide one of the most complicated, among the many constructions of the notion of sign of a permutation.
\end{remark}

\section{The binary case}
\label{binarysec}
In this section, we will focus on $n=2$, i.e., the case of two binary forms $F_1$, $F_2$. We will review the known formulas for the resultant ${\rm Res}_{d_1,d_2}$ and its coefficients, and we will give an example of computation of such coefficients using the method in \S\ref{formulasec}.
We first switch to more convenient notations and write
\[
F_1=a_0 x_1^{d_1}+a_1 x_{1}^{d_1-1} x_2+\cdots+a_{d_1-1} x_1 x_2^{d_1-1}+a_{d_1} x_2^{d_1}\ ,
\]
and
\[
F_2=b_0 x_1^{d_2}+b_1 x_{1}^{d_2-1} x_2+\cdots+b_{d_2-1} x_1 x_2^{d_2-1}+b_{d_2} x_2^{d_2}\ .
\]
It means that we now write $a_i$ for the formal variable earlier denoted by $a_{1,(d_1-i,i)}$, and we write $b_j$ instead of $a_{2,(d_2-j,j)}$.
In this section, we assume knowledge of the standard terminology and notation for integers partitions, in relation to the theory of symmetric functions (see, e.g.,~\cite[Ch. 7]{Stanley}).
As is well known, the resultant of two binary forms is given by the $(d_1+d_2)\times(d_1+d_2)$ determinant
\begin{equation}
{\rm Res}_{d_1,d_2}(F_1,F_2)=\left|
\begin{array}{cccccccc}
a_0 & a_1 & \cdots & \cdots & a_{d_1} & 0 & \cdots & 0 \\
0 & a_0 & a_1 & \cdots & \cdots & a_{d_1} & \ddots & \vdots \\
\vdots & \ddots & \ddots & \ddots & & & \ddots & 0 \\
0 & \cdots & 0 & a_0 & a_1 & \cdots & \cdots & a_{d_1} \\
b_0 & b_1 & \cdots & \cdots & b_{d_2} & 0 & \cdots & 0 \\
0 & b_0 & b_1 & \cdots & \cdots & b_{d_2} & \ddots & \vdots \\
\vdots & \ddots & \ddots & \ddots & & & \ddots & 0 \\
0 & \cdots & 0 & b_0 & b_1 & \cdots & \cdots & b_{d_2}
\end{array}
\right|\ \ ,
\label{Sylvestereq}
\end{equation}
where the $a_0$'s and the $b_{d_2}$'s are aligned along the diagonal.
We will rewrite our expansion of the resultant, in the binary case, as
\[
{\rm Res}_{d_1,d_2}=
\sum_{\substack{\mu,\lambda \\ |\mu|=|\lambda|}} r_{\mu,\lambda}
\ a^{{\rm left}(\mu)}\ b^{{\rm right}(\lambda)}\ ,
\]
where we are summing over pairs of integer partitions $\mu,\lambda$. We impose that $\mu$ is contained in the partition $d_1^{d_2}$, i.e., with $d_2$ parts equal to $d_1$,
and $\lambda$ is contained in the partition $d_2^{d_1}$, i.e., with $d_1$ parts equal to $d_2$. For $1\le i\le d_1$, we denote by $m(\mu,i)$ the multiplicity of the part $i$ in the partition $\mu$. Likewise, for $1\le j\le d_2$, we denote by $m(\lambda,j)$ the multiplicity of the part $j$ in the partition $\lambda$. We use the standard notation $\ell(\cdot)$ for the length of a partition, i.e., the number of nonzero parts, and $|\cdot|$ for the weight of a partition, i.e., the sum of its parts.
We define the multiplicity of the zero part in $\mu$ by $L(\mu):=d_2-\ell(\mu)$.
The same notion for $\lambda$ is slightly different, namely, we
define the multiplicity of the zero part in $\lambda$ by $R(\lambda):=d_1-\ell(\lambda)$.
In the above formula for the resultant we used the notation
\[
a^{{\rm left}(\mu)}:=a_0^{L(\mu)}a_1^{m(\mu,1)}\cdots a_{d_1}^{m(\mu,d_1)}\ ,
\]
and
\[
b^{{\rm right}(\lambda)}:=
b_0^{m(\lambda,d_2)}b_1^{m(\lambda,d_2-1)}\cdots 
b_{d_2-2}^{m(\lambda,2)}
b_{d_2-1}^{m(\lambda,1)} b_{d_2}^{R(\lambda)}\ .
\]
Note that the pairs of partitions $\mu,\lambda$ as above, and which have equal weight, are in bijection with fillings $A$ in ${\rm REF}$. The correspondence is realized, for the first row by
\[
A_{1,(d_1-i,i)}=m(\mu,i)\ ,
\]
for $1\le i\le d_1$, and
\[
A_{1,(d_1,0)}=L(\mu)\ .
\]
For the second row, it is given by
\[
A_{2,(d_2-j,j)}=m(\lambda,d_2-j)\ ,
\]
for $0\le j\le d_2-1$, and
\[
A_{2,(0,d_2)}=R(\lambda)\ .
\]
It is not easy to extract a nice formula for the coeffcients of the resultant from the determinant (\ref{Sylvestereq}). However, with the help of the theory of symmetric functions, one can show (see, e.g.,~\cite[Prop. 16]{GKZAdv}) that
\[
r_{\mu,\lambda}=(-1)^{|\mu|}\ D_{\lambda,\mu}^{-1}\ ,
\]
where $D^{-1}$ is the inverse of the symmetric matrix $D$ indexed by partitions, defined as follows.
We let $D_{\mu,\lambda}$ be the number of $d_2\times d_1$ matrices with entries equal to $0$ or $1$ with row sums given by the parts (including zero parts) of $\mu$, and column sums given by the parts (again including zero parts) of $\lambda$.
By the `Fermionic' variant of the RSK correspondence (see, e.g.,~\cite[Cor. 4.1.15]{Brualdi}), one has
\[
D_{\mu,\lambda}=\sum_{\nu}K_{\rho,\mu} K_{\rho',\lambda}\ ,
\]
where $\rho'$ denotes the transpose of the partition $\rho$, and $K$ is the {\it Kostka matrix}.
Recall that $K_{\rho,\mu}$ by definition counts the number of semistandard Young tableaux of shape $\rho$ and content $\mu$, i.e., such that the number $i$ appears $\mu_i$ times, if $\mu_1\ge\mu_2\ge\cdots$ denote the parts of $\mu$.
As a result, we have the following formula for the coefficients of the binary resultant
\[
r_{\mu,\lambda}=(-1)^{|\mu|}
\sum_{\rho} K_{\mu,\rho}^{-1} K_{\lambda,\rho'}^{-1}\ .
\]
This can be used for concrete computations, via the E\~{g}ecio\~{g}lu-Remmel formula for the inverse Kostka matrix, in terms of special rim hook tableaux~\cite{EgeciogluR}.
The following example with $d_1=6$, $d_2=4$, can be done by hand, with $14$ possibilities for $\rho$, given the dominance order requirement $\mu\succeq\rho\succeq\lambda'$.
The coefficient of the monomial $a_0 a_1 a_3 a_6 b_0^2 b_2 b_4^3$ in ${\rm Res}_{6,4}$ is given by
\begin{eqnarray*}
r_{(6,3,1),(4,4,2)} & = &
  K_{(6,3,1),(6,3,1)}^{-1} K_{(4,4,2),(3,2,2,1,1,1)}^{-1}
+K_{(6,3,1),(6,2,2)}^{-1} K_{(4,4,2),(3,3,1,1,1,1)}^{-1} \\
 & & +K_{(6,3,1),(6,2,1,1)}^{-1} K_{(4,4,2),(4,2,1,1,1,1)}^{-1}
+K_{(6,3,1),(5,4,1)}^{-1} K_{(4,4,2),(3,2,2,2,1)}^{-1} \\
 & & +K_{(6,3,1),(5,3,2)}^{-1} K_{(4,4,2),(3,3,2,1,1)}^{-1}
+K_{(6,3,1),(5,3,1,1)}^{-1} K_{(4,4,2),(4,2,2,1,1)}^{-1} \\
 & & +K_{(6,3,1),(5,2,2,1)}^{-1} K_{(4,4,2),(4,3,1,1,1)}^{-1}
+K_{(6,3,1),(4,4,2)}^{-1} K_{(4,4,2),(3,3,2,2)}^{-1} \\
 & & +K_{(6,3,1),(4,4,1,1)}^{-1} K_{(4,4,2),(4,2,2,2)}^{-1}
+K_{(6,3,1),(4,3,3)}^{-1} K_{(4,4,2),(3,3,3,1)}^{-1} \\
 & & +K_{(6,3,1),(4,3,2,1)}^{-1} K_{(4,4,2),(4,3,2,1)}^{-1}
+K_{(6,3,1),(4,2,2,2)}^{-1} K_{(4,4,2),(4,4,1,1)}^{-1} \\
 & & +K_{(6,3,1),(3,3,3,1)}^{-1} K_{(4,4,2),(4,3,3)}^{-1}
+K_{(6,3,1),(3,3,2,2)}^{-1} K_{(4,4,2),(4,4,2)}^{-1}\ ,
\end{eqnarray*}
which evaluates (without changing the placement of terms and factors) to
\begin{eqnarray*}
r_{(6,3,1),(4,4,2)} & = & 
  1\times 1
+(-1)\times 1 \\
 & & +(-1)\times(-2)
+(-1)\times 0 \\
 & & +0\times(-2)
+0\times 1 \\
 & & +2\times 1
+1\times 1 \\
 & & +1\times(-1)
+(-1)\times 1 \\
 & & +(-1)\times 0
+(-1)\times(-1) \\
 & & +1\times(-1)
+1\times 1\ .
\end{eqnarray*}
Namely,
\[
r_{(6,3,1),(4,4,2)}=4\ .
\]

The method from \S\ref{formulasec} results in quite different computations.
The resolvent equation for a coefficient $r_{\mu,\lambda}$ reads
\begin{eqnarray}
r_{\mu,\lambda} & = & \sum_{\wtmu,\wtla}
\ \ (-1)^{1+L(\mu)+L(\wtmu)+R(\lambda)+R(\wtla)}
\ \ \times\ r_{\wtmu,\wtla} \nonumber \\
 & & \times
\binom{L(\wtmu)}{\ast\ ,m(\mu,2)-m(\wtmu,2),\ldots,m(\mu,d_1)-m(\wtmu,d_1)}
\nonumber \\
 & & \times
\binom{R(\wtla)}{m(\lambda,1)-m(\wtla,1),\ldots,m(\lambda,d_2-1)-m(\wtla,d_2-1),\ast\ }\ .
\label{binresolveq}
\end{eqnarray}
The sum is over pairs of partitions $(\wtmu,\wtla)\neq(\mu,\lambda)$ such that $\wtmu\subset(d_1^{d_2})$,
$\wtla\subset(d_2^{d_1})$, $|\wtmu|=|\wtla|$, together with the inequalities
\begin{eqnarray*}
m(\wtmu,1) & \le & m(\mu,1)+L(\mu) \\
m(\wtmu,2) & \le & m(\mu,2) \\
 & \vdots & \\
m(\wtmu,d_1) & \le & m(\mu,d_1)\ ,
\end{eqnarray*}
and
\begin{eqnarray*}
m(\wtla,1) & \le & m(\lambda,1) \\
 & \vdots & \\
m(\wtla,d_2-1) & \le & m(\lambda,d_2-1) \\
m(\wtla,d_2) & \le & m(\lambda,d_2)+R(\lambda)\ .
\end{eqnarray*}
In the two multinomial coefficients appearing in (\ref{binresolveq}), the symbol $\ast$ stands for the value which makes the sum of bottom entries equal to the top entry.

The computation of the same coefficient $r_{(6,3,1),(4,4,2)}$ now requires the values of $6$ other coefficients, including $r_{\varnothing,\varnothing}=1$, due to the initial condition (\ref{initcondeq}). The notation $\varnothing$ is for the empty partition only made of zero parts. Note that $r_{\varnothing,\varnothing}$ corresponds to the monomial $a_0^{d_2} b_{d_2}^{d_1}$, i.e., the contribution of the diagonal to the determinant (\ref{Sylvestereq}).  
The needed resolvent equations, and the initial condition, form the following system
\[
\begin{pmatrix}
1 & -3 & -2 & -4 & -10 & 0 & 72 \\
0 & 1 & 2 & 4 & 0 & 0 & -24 \\
0 & 0 & 1 & 0 & 0 & 5 & -6 \\
0 & 0 & 0 & 1 & 0 & 0 & -4 \\
0 & 0 & 0 & 0 & 1 & 0 & -4 \\
0 & 0 & 0 & 0 & 0 & 1 & -1 \\
0 & 0 & 0 & 0 & 0 & 0 & 1 
\end{pmatrix}
\begin{pmatrix}
r_{(6,3,1),(4,4,2)} \\
r_{(6),(4,2)} \\
r_{(1,1),(2)} \\
r_{(6,1,1),(4,4)} \\
r_{(3,1),(4)} \\
r_{(1,1,1,1),(4)} \\
r_{\varnothing,\varnothing}
\end{pmatrix}
=
\begin{pmatrix}
0 \\
0 \\
0 \\
0 \\
0 \\
0 \\
1
\end{pmatrix}\ .
\]
Its solution is
\[
\begin{pmatrix}
r_{(6,3,1),(4,4,2)} \\
r_{(6),(4,2)} \\
r_{(1,1),(2)} \\
r_{(6,1,1),(4,4)} \\
r_{(3,1),(4)} \\
r_{(1,1,1,1),(4)} \\
r_{\varnothing,\varnothing}
\end{pmatrix}
=
\begin{pmatrix}
4 \\
6 \\
1 \\
4 \\
4 \\
1 \\
1
\end{pmatrix}\ ,
\]
from which we again find $
r_{(6,3,1),(4,4,2)}=4$.

The above two computations, via special rim hook tableaux, and via \S\ref{formulasec}, were done by hand, and they seemed comparable in the amount of effort.
It would be interesting to do a more systematic comparison of the two algorithms, from the point of view of complexity theory.

\begin{remark}
Given a linear system $Ax=y$, involving a square invertible matrix $A$, showing the necessary condition $x=A^{-1}y$ ($x$ must be given by $A^{-1}y$) involves the identity $A^{-1}A={\rm I}$, whereas showing the sufficient condition $x=A^{-1}y$ ($A^{-1}y$ solves the system) involves the identity $A A^{-1}={\rm I}$.
In this article, leveraging the available theory of resultants, we showed the necessary condition that the $r_A$ must be given by the formula in Theorem \ref{maintheorem}.
We are at present unable to directly prove the sufficient condition that the $r_A$ defined by the formula in Theorem \ref{maintheorem} are coefficients of an IFDZ of appropriate multidegree. This would require understanding syzygies among the equations ${\rm Eq}(\mathbb{E},C,\beta)$.
In this regard, it is noteworthy that the identity $KK^{-1}={\rm I}$, for the Kostka matrix, can be proved using the special rim hook tableau formula~\cite[\S3]{EgeciogluR}, whereas a direct combinatorial proof of $K^{-1}K={\rm I}$ is still missing, as far as we know.
\end{remark}

\section{Bound on the height of resultants}
\label{heightsec}

In this section we show how one can use our main result in order to
recover the following bound on the height of resultants, which is well known (see~\cite{Sombra}), and as far as we are aware, is the best available general bound which holds for arbitrary dimension $n$ and degrees $d_1,\ldots,d_n$.
Before stating the bound, recall that if $G$ is a polynomial with complex coefficients, the height of $G$ denoted by $H(G)$ is the maximum modulus of a coefficients of $G$.
We will also denote by $L(G)$ the sum of the moduli of the coefficients. Namely, $H(G)$ is the $\ell^{\infty}$ norm, while $L(G)$ is the $\ell^1$ norm of the coefficients seen as complex-valued functions of multiindices encoding the monomials.

\begin{theorem}\label{heightthm}
For any $n\ge 1$, and any degrees $d_1,\ldots,d_n$, we have
\[
H({\rm Res}_{d_1,\ldots,d_n}) \le\prod_{i=1}^{n}
{\binom{d_i+n-1}{n-1}}^{\delta_i}\ ,
\]
where, like earlier in this article, the deltas are defined by $\delta_i=\prod_{\substack{j=1\\
j\neq i}}^{n}d_j$.
\end{theorem}

We will give a new proof of this theorem, by induction on $n$. We do this by first showing what we call an admissible bound
\[
H({\rm Res}_{d_1,\ldots,d_n}) \le \mathcal{B}_{d_1,\ldots,d_n}\ ,
\]
with strictly positive right-hand side, 
for fixed $n$, but for general degrees $d_1,\ldots,d_n$.
By admissible bound we mean one which has the property
\[
\lim\limits_{k\rightarrow\infty}\frac{1}{k^n}\ln \mathcal{B}_{kd_1,\ldots,kd_n}=0\ .
\]
By a well known argument, e.g., used by Sombra in~\cite{Sombra} (see also~\cite{GalaJetal}), once an admissible bound is established, the inequality in the above theorem follows. By Fourier inversion for multiple Fourier series, one easily sees that the height $H(G)$ of a polynomial is bounded by the supremum of evaluations $|G(u)|$ where the variables are specialized to complex numbers of unit modulus. Moreover, if $G(u_1,\ldots,u_N)$ is multihomogeneous in $N$ series of variables $u_i=(u_{i,j})_{1\le j\le k_i}$, $1\le i\le N$, of respective degrees $D_i$, then one easily obtains the elementary bound
\begin{equation}
|G(u_1,\ldots,u_N)|\le H(G)\times
\prod_{i=1}^{N}
\left(
\sum_{j=1}^{k_i}|u_{i,j}|
\right)^{D_i}\ ,
\label{Guboundeq}
\end{equation}
where the $u$ variables have been specialized to arbitrary complex numbers.
By the multiplicative property of resultants, the bound (\ref{Guboundeq})
applied to $G={\rm Res}_{kd_1,\ldots,kd_n}$, as well as the easy inequality $L(F_i^k)\le L(F_i)^k$,
one readily obtains
\begin{eqnarray*}
\left|
{\rm Res}_{d_1,\ldots,d_n}(F_1,\ldots,F_n)
\right|^{k^n}&=&\left|
{\rm Res}_{kd_1,\ldots,kd_n}(F_1^k,\ldots,F_n^k)
\right| \\
 & \le &
H({\rm Res}_{kd_1,\ldots,kd_n}) \times\prod_{i=1}^{n}
{\binom{d_i+n-1}{n-1}}^{k^n\delta_i}\ .
\end{eqnarray*}
Here, the forms $F_i$ are specialized to ones with coefficients given by complex numbers of unit modulus.
We then have
\[
\left|
{\rm Res}_{d_1,\ldots,d_n}(F_1,\ldots,F_n)
\right|\le
\left(\mathcal{B}_{kd_1,\ldots,kd_n}\right)^{\frac{1}{k^n}}
\times
\prod_{i=1}^{n}
{\binom{d_i+n-1}{n-1}}^{\delta_i}\ .
\]
Then, taking the supremum over forms $F_i$ whose coefficients are complex numbers of unit modulus, we deduce
\[
H({\rm Res}_{d_1,\ldots,d_n})
\le
\left(\mathcal{B}_{kd_1,\ldots,kd_n}\right)^{\frac{1}{k^n}}
\times
\prod_{i=1}^{n}
{\binom{d_i+n-1}{n-1}}^{\delta_i}\ .
\]
Hence, by taking the $k\rightarrow\infty$ limit, we see that admissibility in the seed bound with $\mathcal{B}$ allows one to deduce the improved bound in Theorem~\ref{heightthm}.

We now proceed to establishing the required admissible bound. For better readability we will write $A(i,\alpha)$ instead of $A_{i,\alpha}$ for entries of fillings.
We also introduce the notation $\mathbb{S}:=\mathbb{D}\backslash(\mathbb{E}\cup\mathbb{T})$, and $\mathbb{S}_i:=\widetilde{\mathbb{D}}_i\backslash\{\varepsilon_{i\ast},\tau_{i\ast}\}$. Namely, $\mathbb{S}_i$ is the set of multiindices $\alpha\in\mathbb{N}^n$, with $|\alpha|=d_i$, other than $\varepsilon_{i\ast}$ and $\tau_{i\ast}$ coming from the transversals $\mathbb{E}$ and $\mathbb{T}$.
Going back to the transition matrix defined in (\ref{Tdefeq}), we show that it satisfies the following bound.

\begin{lemma}\label{basicTlem}
For $A\in{\rm REF}\backslash{\rm FRR}$ and $\wtA\in{\rm FRR}$, we have
\[
|\mathscr{T}(A,\wtA)|\le
\bbone\{\mathbf{1}_{\mathbb{S}}\wtA\le\mathbf{1}_{\mathbb{S}}A\}\times
\frac{\prod_{i=1}^{n}\left(\delta_i-\sum_{\alpha\in\mathbb{S}_i}\wtA(i,\alpha)\right)!}{\prod_{i=1}^{n}\left(\delta_i-\sum_{\alpha\in\mathbb{S}_i}A(i,\alpha)\right)!\times
\prod_{(i,\alpha)\in\mathbb{S}}\left(A(i,\alpha)-\wtA(i,\alpha)\right)!}\ .
\]
\end{lemma}

\noindent{\bf Proof:}
We first note that the $i$-th component of
${\rm Row}(\mathbf{1}_{\mathbb{E}}\wtA)$ is
\[
\wtA(i,\varepsilon_{i\ast})=\delta_i-\wtA(i,\tau_{i\ast})-\sum_{\alpha\in\mathbb{S}_i}\wtA(i,\alpha)\ ,
\]
because of the imposed homogeneity of degree $\delta_i$, with respect to the coefficients of the form $F_i$.
The entry at a cell $(i,\alpha)\in\mathbb{D}$ of the filling $\mathbf{1}_{\mathbb{D}\backslash\mathbb{E}}(A-\wtA)+P_{\mathbb{T}}(\mathbf{1}_{\mathbb{E}}A)$
is equal to zero if $(i,\alpha)\in\mathbb{E}$. It is equal to $A(i,\alpha)-\wtA(i,\alpha)$, if $(i,\alpha)\in\mathbb{S}$. Finally, if $(i,\alpha)\in\mathbb{T}$, it is equal to
\[
A(i,\tau_{i\ast})+A(i,\varepsilon_{i\ast})-\wtA(i,\tau_{i\ast})=
\delta_i-\wtA(i,\tau_{i\ast})-\sum_{\alpha\in\mathbb{S}_i}A(i,\alpha)\ .
\]
Therefore,
\begin{eqnarray*}
\lefteqn{
\begin{bmatrix}
{\rm Row}(\mathbf{1}_{\mathbb{E}}\wtA) \\
\mathbf{1}_{\mathbb{D}\backslash\mathbb{E}}(A-\wtA)+P_{\mathbb{T}}(\mathbf{1}_{\mathbb{E}}A)
\end{bmatrix}
=} & & \\
 & & \frac{\prod_{i=1}^{n}\left(\delta_i-\wtA(i,\tau_{i\ast})
-\sum_{\alpha\in\mathbb{S}_i}\wtA(i,\alpha)\right)!}{\prod_{i=1}^{n}\left(\delta_i
-\wtA(i,\tau_{i\ast})
-\sum_{\alpha\in\mathbb{S}_i}A(i,\alpha)\right)!\times
\prod_{(i,\alpha)\in\mathbb{S}}\left(A(i,\alpha)-\wtA(i,\alpha)\right)!}\ .
\end{eqnarray*}
Recall that binomial coefficients $\binom{m}{p}$ increase with $m$ when $p$ is fixed.
As a result,
\[
\frac{(a-u)!}{(b-u)!}=(a-b)!\binom{a-u}{a-b}
\]
decreases with $u$, when $a\ge b$. Bounding the latter ratio by the value at $u=0$, we see that
\[
\frac{(a-u)!}{(b-u)!}\le\frac{a!}{b!}\ ,
\]
when $a\ge b\ge u\ge 0$.
For every $i\in[n]$, we use the above inequality for
\begin{eqnarray*}
a &= & \delta_i-\sum_{\alpha\in\mathbb{S}_i}\wtA(i,\alpha) \\
b & = & \delta_i-\sum_{\alpha\in\mathbb{S}_i}A(i,\alpha) \\
u & = & \wtA(i,\tau_{i\ast})\ .
\end{eqnarray*}
We also partially relax the condition in the indicator function, by only taking cells in $\mathbb{S}$ into consideration, and the lemma follows. 
\qed

\begin{lemma}\label{Tpowerlem}
For $k\ge 0$, $A\in{\rm REF}\backslash{\rm FRR}$ and $B\in{\rm FRR}$, we have
\begin{eqnarray*}
|\mathscr{T}^{k+1}(A,B)| & \le &
\bbone\{\mathbf{1}_{\mathbb{S}}\wtA\le\mathbf{1}_{\mathbb{S}}A\}
\times\prod_{i=1}^{n}(\delta_i+1)^{k} \\
 & & \times
\frac{\prod_{i=1}^{n}\left(\delta_i-\sum_{\alpha\in\mathbb{S}_i}B(i,\alpha)\right)!\times
\prod_{(i,\alpha)\in\mathbb{S}}(k+1)^{A(i,\alpha)-B(i,\alpha)}}{\prod_{i=1}^{n}\left(\delta_i-\sum_{\alpha\in\mathbb{S}_i}A(i,\alpha)\right)!\times
\prod_{(i,\alpha)\in\mathbb{S}}\left(A(i,\alpha)-B(i,\alpha)\right)!}\ .
\end{eqnarray*}
\end{lemma}

\noindent{\bf Proof:}
We proceed by induction on $k$, with the $k=0$ case already taken care of by Lemma \ref{basicTlem}, so we now assume $k\ge 1$. From the definition (\ref{Tpowereq}), we see that
\[
\mathscr{T}^{k+1}(A,B)=\sum_{\wtA\in{\rm REF}\backslash{\rm FRR}}
\bbone\{B\prec\wtA\prec A\}
\mathscr{T}(A,\wtA)\ \mathscr{T}^{k}(\wtA,B)\ .
\]
As a result of Lemma \ref{basicTlem} and the induction hypothesis, we obtain
\begin{eqnarray*}
\left|\mathscr{T}^{k+1}(A,B)\right| & \le & 
\sum_{\wtA\in{\rm EF}(\delta)}
\bbone\{\mathbf{1}_{\mathbb{S}}B\le
\mathbf{1}_{\mathbb{S}}\wtA\le\mathbf{1}_{\mathbb{S}}A\} \\ 
 & & \times
\frac{\prod_{i=1}^{n}\left(\delta_i-\sum_{\alpha\in\mathbb{S}_i}\wtA(i,\alpha)\right)!}{\prod_{i=1}^{n}\left(\delta_i-\sum_{\alpha\in\mathbb{S}_i}A(i,\alpha)\right)!\times
\prod_{(i,\alpha)\in\mathbb{S}}\left(A(i,\alpha)-\wtA(i,\alpha)\right)!} \\
 & & \times
\frac{\prod_{i=1}^{n}(\delta_i+1)^{k-1}\times
\prod_{i=1}^{n}\left(\delta_i-\sum_{\alpha\in\mathbb{S}_i}B(i,\alpha)\right)!\times
\prod_{(i,\alpha)\in\mathbb{S}} k^{\wtA(i,\alpha)-B(i,\alpha)}}{\prod_{i=1}^{n}\left(\delta_i-\sum_{\alpha\in\mathbb{S}_i}\wtA(i,\alpha)\right)!\times
\prod_{(i,\alpha)\in\mathbb{S}}\left(\wtA(i,\alpha)-B(i,\alpha)\right)!}\ .
\end{eqnarray*}
Note that we also enlarged the summation set to ${\rm EF}(\delta)$ which, in particular, means we now ignore the weight condition.
For each cell $(i,\alpha)\in\mathbb{S}$, we must sum over the integer $\wtA(i,\alpha)$ 
over the range from $B(i,\alpha)$ to $A(i,\alpha)$. 
We do this using the Newton binomial theorem which gives
\[
\sum_{\wtA(i,\alpha)=B(i,\alpha)}^{A(i,\alpha)}
\frac{k^{\wtA(i,\alpha)-B(i,\alpha)}}{\left(A(i,\alpha)-\wtA(i,\alpha)\right)!\ 
\left(\wtA(i,\alpha)-B(i,\alpha)\right)!}
=\frac{(k+1)^{A(i,\alpha)-B(i,\alpha)}}{\left(A(i,\alpha)-B(i,\alpha)\right)!}\ .
\]
For each $i\in[n]$, we must also sum over the integer $\wtA(i,\tau_{i\ast})$
for which we allow the largest possible range from $0$ to $\delta_i$. This produces an additional factor $\prod_{i=1}^{n}(\delta_i+1)$.
A little algebra and cleaning up the resulting upper bound gives the desired inequality for $k$, and the lemma follows by induction.
\qed

We now use the previous lemma to get an estimate on the quantity $\mathscr{C}(A,B)$
defined in (\ref{CABdefeq}).
Recall that the sum over $k$ terminates because of the limitation
\[
k\le \mathcal{Z}(A)=\sum_{(i,\alpha)\in\mathbb{S}} A(i,\alpha)\le \sum_{i=1}^{n} \delta_i\ .
\]
From the multinomial theorem, and bounding a single term by an entire sum of nonnegative terms, we have the trivial inequality
\[
\binom{v_i}{u_1,\ldots,u_{m_i}}\le m_{i}^{v_i}\ .
\]
We use this, for each $i\in[n]$, with
\[
v_i=\delta_i-\sum_{(i,\alpha)\in\mathbb{S}}B(i,\alpha)\le \delta_i\ ,
\]
and
\[
m_i=|\mathbb{S}_i|\le|\widetilde{D}_i|=\binom{d_i+n-1}{n-1}\ .
\]
This allows us to bound the big ratio of factorials in Lemma \ref{Tpowerlem} by
\[
\frac{\prod_{i=1}^{n}\left(\delta_i-\sum_{\alpha\in\mathbb{S}_i}B(i,\alpha)\right)!}{\prod_{i=1}^{n}\left(\delta_i-\sum_{\alpha\in\mathbb{S}_i}A(i,\alpha)\right)!\times
\prod_{(i,\alpha)\in\mathbb{S}}\left(A(i,\alpha)-B(i,\alpha)\right)!}
\le \prod_{i=1}^{n}
{\binom{d_i+n-1}{n-1}}^{\delta_i}\ .
\]
We also use the bound
\[
\prod_{(i,\alpha)\in\mathbb{S}}(k+1)^{A(i,\alpha)-B(i,\alpha)}
\le (\delta_1+\cdots+\delta_n+1)^{\sum_{(i,\alpha)\in\mathbb{S}} A(i,\alpha)}
\le (\delta_1+\cdots+\delta_n+1)^{\delta_1+\cdots+\delta_n}\ .
\]
The sum over $k$ then gives
\[
\sum_{k=0}^{\mathcal{Z}(A)}\prod_{i=1}^{n}(\delta_i+1)^k
=\frac{\prod_{i=1}^{n}(\delta_i+1)^{\mathcal{Z}(A)+1}-1}{\prod_{i=1}^{n}(\delta_i+1)-1}
\le \prod_{i=1}^{n}(\delta_i+1)^{\delta_1+\cdots+\delta_n+1}\ .
\]
We have thus showed that, for all $A\in{\rm REF}$, and all $B\in{\rm FRR}$, we have
\begin{eqnarray}
|\mathscr{C}(A,B)| & \le &
2\times\bbone\{\mathbf{1}_{\mathbb{S}}B\le \mathbf{1}_{\mathbb{S}}A\}
\times \prod_{i=1}^{n}
{\binom{d_i+n-1}{n-1}}^{\delta_i} \nonumber \\
 & & \times (\delta_1+\cdots+\delta_n+1)^{\delta_1+\cdots+\delta_n}
\times \prod_{i=1}^{n}(\delta_i+1)^{\delta_1+\cdots+\delta_n+1}\ .
\label{CABboundeq}
\end{eqnarray}
The first factor of $2$ is due to the $A=B$ term which is much smaller than our bound on the $\sum_{k\ge 0}$ contribution to (\ref{CABdefeq}).
From formula (\ref{FRRformulaeq}), we have that for all $B\in{\rm FRR}$,
\[
|r_B|\le
H({\rm Res}_{d_2,\ldots,d_n})^{d_1}\times
|{\rm REF}^{(2)}|^{d_1} \ .
\]
We use the crude bound
\[
|{\rm REF}^{(2)}|\le\prod_{i=2}^{n}
\binom{\delta_i^{(2)}+|\widetilde{D}_{i}^{(2)}|-1}{|\widetilde{D}_{i}^{(2)}|-1}
\le \prod_{i=2}^{n} 2^{\delta_i^{(2)}+|\widetilde{D}_{i}^{(2)}|-1}
\]
which is only due to the multihomogenity of the $(n-1)$-dimensional resultant ${\rm Res}_{d_2,\ldots,d_n}$.
We also have $|\widetilde{D}_{i}^{(2)}|\le (d_i+1)^{n-2}$, since we are counting multiindices $(0,\alpha_2,\ldots,\alpha_n)$ of length $d_i$, which are determined by $\alpha_2,\ldots,\alpha_{n-1}$ which range over $\{0,1,\ldots,d_i\}$.
As a result, we get for any $B\in{\rm FRR}$,
\begin{eqnarray}
|r_B| & \le & 
\left[
\prod_{i=2}^{n}
{\binom{d_i+n-2}{n-2}}^{\delta_i^{(2)}}
\right]^{d_1}
\times \left[
\prod_{i=2}^{n}
2^{\delta_i^{(2)}+(d_i+1)^{n-2}-1}
\right]^{d_1} \nonumber \\
 & \le & \prod_{i=2}^{n}
{\binom{d_i+n-2}{n-2}}^{\delta_i}
\times
\prod_{i=2}^{n}
2^{\delta_i+d_1 (d_i+1)^{n-2}}\ ,
\label{rBboundeq}
\end{eqnarray}
after noting that $\delta_i=d_1\delta_{i}^{(2)}$, for $2\le i\le n$, and using the induction hypothesis, on the dimension $n$, for our proof of Theorem \ref{heightthm}.
We finally use
\[
\binom{d_i+n-2}{n-2}\le
\binom{d_i+n-1}{n-1}
\le (d_i+1)^{n-1}\ ,
\]
and combine (\ref{CABboundeq}) and (\ref{rBboundeq}), with a slight increase in the range of the product of binomials in the latter,
in order to arrive at
\[
H({\rm Res}_{d_1,\ldots,d_n})\le\mathcal{B}_{d_1,\ldots,d_n}\ ,
\]
where
\begin{eqnarray*}
\mathcal{B}_{d_1,\ldots,d_n} &:= & 2\times (\delta_1+\cdots+\delta_n+1)^{\delta_1+\cdots+\delta_n}
\times \prod_{i=1}^{n}(\delta_i+1)^{\delta_1+\cdots+\delta_n+1} \\
 & & \times \prod_{i=1}^{n}(d_i+1)^{2(n-1)\delta_i}
\times\prod_{i=2}^{n} 2^{\delta_i+d_1 (d_i+1)^{n-2}}\ .
\end{eqnarray*}
It is easy to see that this bound is admissible, because the exponents are of degree $n-1$ in the form degrees $d_1,\ldots,d_n$. Namely, when these degrees are scaled by a factor of $k$, the amplification ratio $k^{n-1}$ gets crushed once we divide by $k^n$ after taking the logarithm.
This concludes the proof of Theorem \ref{heightthm}.

\section{Outlook}
\label{outlooksec}

As far as we know, Theorem \ref{maintheorem} provides the first explicit formula for the coefficients of resultants in arbitrary dimension and for general degrees. 
However, the formula which is an alternating sum of products of multinomial coefficients is, unsurprisingly, very complicated. Moreover, it is not yet clear what the multinomials are counting. We believe the most pressing task ahead is to define suitable combinatorial structures counted by the products of multinomials appearing in our formula. Ideally, one would like an expression similar to the rim hook tableau expansion for the inverse Kostka matrix due to E\~{g}ecio\~{g}lu and Remmel~\cite{EgeciogluR}. Namely, we would like a formula 
\begin{equation}
r_A=\sum_{c\in C}{\rm sgn}(c)
\label{altsumformula}
\end{equation}
for the coefficients of resultants, where $C$ is a set of combinatorial objects and ${\rm sgn}(c)=\pm 1$ is a suitable sign assigned to each instance $c$ of such combinatorial structures. We believe this should be doable by carefully going through the steps of the proof of Theorem \ref{maintheorem}. One should also first focus on the binary case $n=2$, and even more modestly, consider the situation where one of the degrees say $d_1$ is small (e.g., equal to 2 or 3) while $d_2$ can be arbitrarily large. Such an investigation could benefit from the related study done in~\cite{DAndreaH}.

Note that a formula (\ref{altsumformula}) may open the door to exploiting sign cancellations, e.g., via sign-reversing involutions, with the aim of obtaining better bounds for the coefficients of resultants. The study of improved bound on heights of resultants is an active area of research (see, e.g.,~\cite{Sombra,DAndreaKS,FernandezP}), and it could benefit from a better combinatorial understanding of the coefficients $r_A$ and the sign cancellations in explicit formulas for the latter. In \S\ref{heightsec}, we were able to recover some of the best available bounds of heights of resultants, without exploiting sign cancellations,
since in Lemma \ref{basicTlem} we threw away the sign from (\ref{Tdefeq}). There is thus much room left for improvement. Note that obtaining what we called an admissible bound $\mathcal{B}_{d_1,\ldots,d_n}$ is not trivial. For example, the height bound resulting from~\cite[Eq. 4.3.41]{BostGS} is not admissible because of the term involving the Stoll number $\sigma_N$.

Another interesting continuation of the present work is to explore the generalization to multigraded resultants where the variables $x_1,\ldots,x_n$ are partitioned into several groups, and the forms $F_i$ are separately homogeneous with respect to each of these groups of variables. Such multigraded resultants have been considered classically by Sylvester, Muir and Lasker (see~\cite{McCoy}). For more recent work see~\cite{SturmfelsZ} and~\cite{BuseCN}.
We are confident that one should be able to generalize the explicit formula in Theorem \ref{maintheorem} to the multigraded setting. This is because 
of the availability of a theory of inertia forms in the multigraded setting~\cite{AwaneCG,BuseCN}.


\end{document}